\documentclass{amsart}
\usepackage{amsmath,amsfonts,mathrsfs,amssymb,enumerate,graphicx,subfig,caption}
\usepackage{xcolor} 
\usepackage[backref=page]{hyperref} 
\renewcommand*{\backref}[1]{}
\renewcommand*{\backrefalt}[4]{%
    \ifcase #1 (Not cited.)%
    \or        (Cited on page~#2.)%
    \else      (Cited on pages~#2.)%
    \fi}

\numberwithin{equation}{section}

\newtheorem{theorem}{Theorem}

\newtheorem{question}{Question}

\newtheorem{proposition}{Proposition}[section]
\newtheorem{corollary}[proposition]{Corollary}
\newtheorem{lemma}[proposition]{Lemma}

\newtheorem{definition}[proposition]{Definition}
\theoremstyle{remark}

\DeclareMathOperator{\im}{im}

\DeclareMathOperator{\rank}{rank}

\DeclareMathOperator{\dimaff}{dim_{\mathsf{aff}}}
\DeclareMathOperator{\dimaffQ}{dim_{\mathsf{aff}}^Q}

\DeclareMathOperator{\dimh}{dim_{\mathsf{H}}}

\DeclareMathOperator{\GL}{GL}
\DeclareMathOperator{\SO}{SO}

\renewcommand{\O}{\mathrm{O}}

\DeclareMathOperator{\diam}{diam}

\newcommand{\threebar}[1]{{\left\vert\kern-0.25ex\left\vert\kern-0.25ex\left\vert #1 
    \right\vert\kern-0.25ex\right\vert\kern-0.25ex\right\vert}}
\newcommand{\iii}{\mathtt{i}}
\newcommand{\jjj}{\mathtt{j}}

\newcommand{\R}{\mathbb{R}}

\newcommand{\A}{\mathsf{A}}
\newcommand{\B}{\mathsf{B}}

\newcommand{\cagri}{\c{C}a\u{g}r\i{ }}

\usepackage{soul}

\begin{document}

\title{Exceptional projections of self-affine sets: an introduction}
\author{Ian D. Morris}
\address{School of Mathematical Sciences, Queen Mary University of London, Mile End Road, London E1 4NS, U.K.}
\email{i.morris@qmul.ac.uk }
\begin{abstract}
We describe some recent results on the dimensions of linear projections of self-affine fractals, focusing in particular on an upper bound for the dimension of the projected image. We give a self-contained treatment of this bound and illustrate it through explicit examples, in the process exhibiting some smooth submanifolds of the Grassmannian which can be contained in the exceptional set in Marstrand's theorem.
\end{abstract}
\maketitle
\section{Introduction and background}
\subsection{Hausdorff dimension and iterated function systems}  \label{ss:onepointone}

If $X$ is a subset of $\R^d$ then the \emph{Hausdorff dimension} of $X$, which we will denote by $\dimh X$, is defined to be the infimal value of $s>0$ such that the quantity
\[\mathcal{H}_\delta^s:=\inf\left\{\sum_{i=1}^\infty (\diam U_i)^s \colon X\subseteq \bigcup_{i=1}^\infty U_i\text{ and }\sup_{i} \diam U_i\leq \delta\right\}\]
is bounded uniformly with respect to $\delta>0$. It is an easy exercise to show that the Hausdorff dimension of familiar and simple sets --- such as a nonempty open subset of $\R^d$, for example --- coincides with the value suggested by one's everyday intuition. A well-known source of sets with non-integer dimension, on the other hand, is the following celebrated result of J.E. Hutchinson:
\begin{theorem}\label{th:hutch}
Let $T_1,\ldots,T_N \colon \R^d \to \R^d$ be contractions of $\R^d$. Then there exists a unique nonempty compact set $X\subset \R^d$ such that $X=\bigcup_{i=1}^N X$. If additionally every $T_i$ is a similarity transformation which contracts the Euclidean norm by precisely $r_i \in (0,1)$, and if there exists a nonempty open set $U\subset \R^d$ such that $T_1U,\ldots,T_NU$ are pairwise disjoint subsets of $U$, then the dimension $s:=\dimh X$ satisfies the equation $\sum_{i=1}^N r_i^s=1$.
\end{theorem}
To avoid trivialities, here and in the sequel we will always assume that $N \geq 2$. By saying that $T_1,\ldots,T_N$ are contractions of $\R^d$ we mean that for some norm $\threebar{\cdot}$ on $\R^d$ and some real number $\lambda \in (0,1)$ one has $\threebar{T_ix-T_iy} \leq \lambda \threebar{x-y}$ for every $i=1,\ldots,N$ and $x,y \in \R^d$. We will reserve the symbol $\|\cdot\|$ to refer specifically to the Euclidean norm on $\R^d$ and to its induced operator norm on the vector space $M_d(\R)$ of all $d \times d$ real matrices. To simplify what follows, given $N \geq 2$ we will say that a \emph{word} is a finite sequence $\iii=i_1 i_2 \cdots i_n$ where each $i_j \in \{1,\ldots,N\}$, and we denote the length $n$ of such a word by $|\iii|$. If $T_1,\ldots,T_N$ are contractions of $\R^d$ and $\iii=i_1\cdots i_n$ is a word then we define $T_\iii$ to be the composition $T_{i_1}\circ T_{i_2}\circ \cdots\circ T_{i_n}$, and similarly if $r_1,\ldots,r_N$ are as above we let $r_\iii$ denote the product $r_{i_1}\cdots r_{i_n}$.

The tuple $(T_1,\ldots,T_N)$ in Theorem \ref{th:hutch} is conventionally called an \emph{iterated function system} and the set $X$ its \emph{attractor}. The attractors of iterated function systems in the case where every $T_i$ is a similarity transformation are generically referred to as \emph{self-similar sets}, and this class of sets encompasses many familiar constructions of fractal sets such as the middle-third Cantor set, Sierpi\'nski triangle, von Koch curve and Menger sponge. In the case where the clause concerning the open set $U$ is not satisfied, the conclusion $\sum_{i=1}^N r_i^s=1$ is in general false, as can easily be seen by considering the case in which all of the transformations $T_1,\ldots,T_N$ are identical: in this case the attractor is a point (hence zero-dimensional) but the solution to $\sum_{i=1}^N r_i^s=1$ is nonzero. On the other hand when no equations of the form $T_\iii =T_\jjj$ occur with $\iii, \jjj$ distinct, when the transformations $T_1,\ldots,T_N$ do not preserve a proper affine subspace of $\R^d$, and when the unique solution $s$ to Hutchinson's equation $\sum_{i=1}^N r_i^s=1$ does not exceed $d$, it is very often the case that the dimension of the attractor solves $\sum_{i=1}^N r_i^s=1$: see for example \cite{Ho14,Ho20}.

When the transformations $T_i$ are allowed to be non-conformal, the problem of determining the dimension of the attractor becomes far more difficult: this direction of extension of Hutchinson's theorem is a major topic of contemporary research (see for example \cite{BaHoRa19,BaKe24,DaSi17,Fe23,FrJu21,MoSe23,Ra24}). By far the most substantially investigated non-conformal case is that in which the transformations $T_i$ are \emph{affine}: that is, for every $i=1,\ldots,N$ there exist $A_i \in M_d(\R)$ and $v_i \in \R^d$ such that $T_ix = A_ix+v_i$ for all $x \in \R^d$. In this work we will always make the stronger assumption that every $A_i$ is additionally invertible, i.e. that $A_i \in \GL_d(\R)$ for every $i$. For such iterated function systems, given a word $\iii=i_1\cdots i_n$ we will write $A_\iii$ to denote the composition $A_{i_1}\cdots A_{i_n}$.

\begin{figure}
    \centering
    \subfloat[The standard Sierpi\'nski triangle covered by a single ball $\mathbf{B}$. ]{{\includegraphics[width=5.7cm]{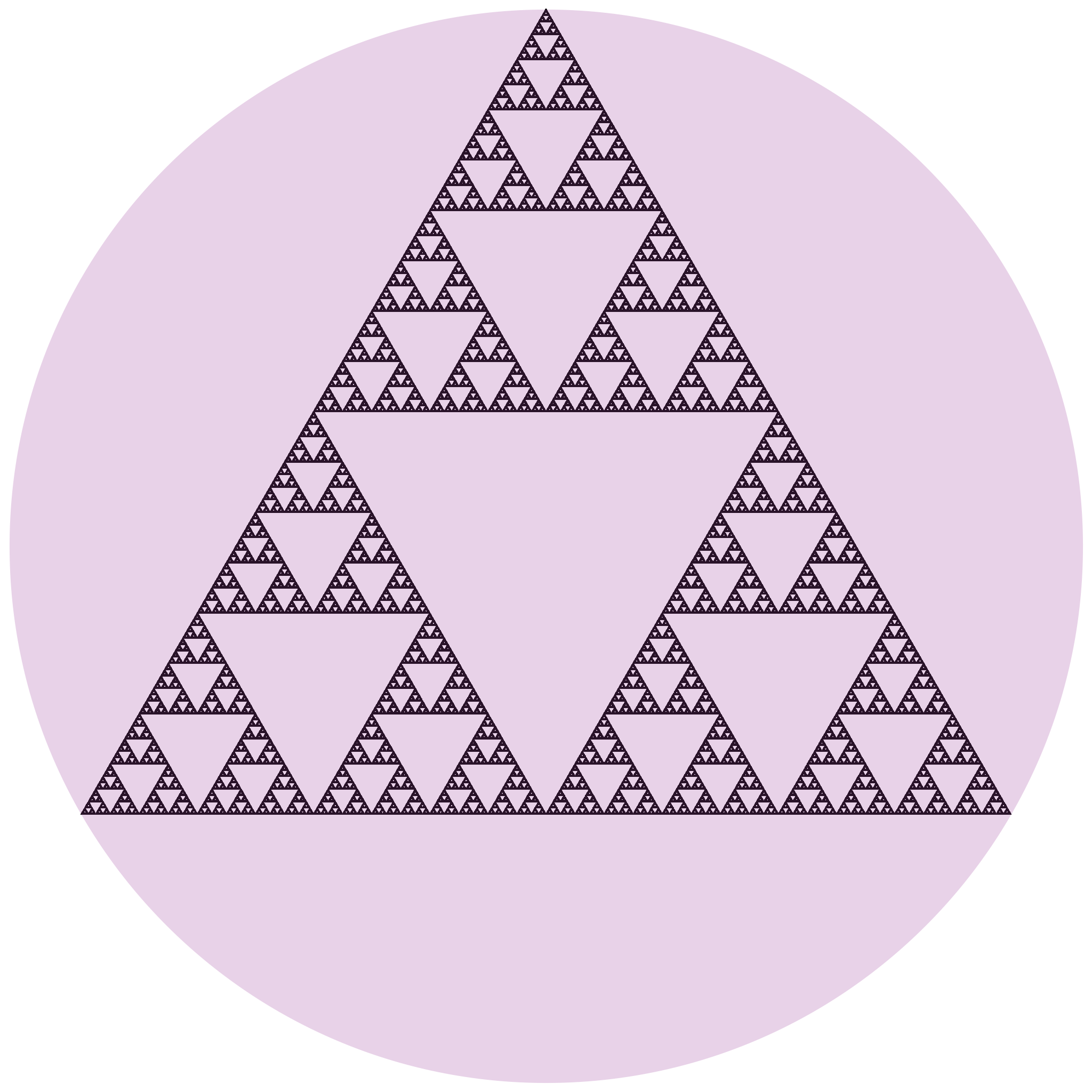}}}
    \quad
    \subfloat[The Sierpi\'nski triangle covered by the images $T_\iii \mathbf{B}$ where $|\iii|=2$.]{{\includegraphics[width=5.7cm]{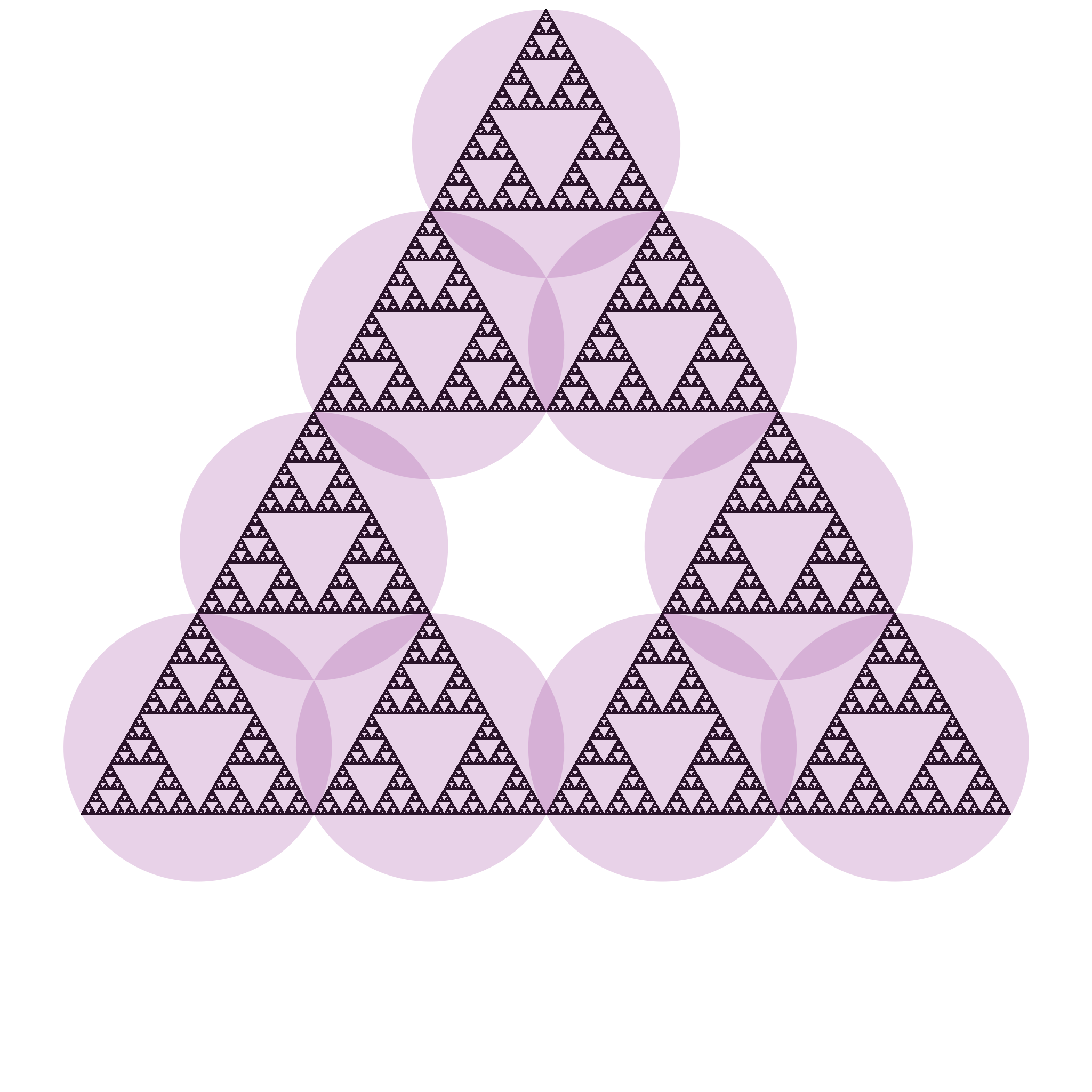}}}
    \caption{Self-similar sets satisfying the hypotheses of Theorem \ref{th:hutch} may be efficiently covered using sets of the form $T_\iii \mathbf{B}$ where $\mathbf{B}$ is a ball which contains the attractor. }
    \label{fi:covers}%
\end{figure}

 The greater difficulty of non-conformal affine examples is at an intuitive level not difficult to explain. If every $T_i$ is a similarity transformation, then given a covering of the attractor $X$ by some collection of balls $\mathbf{B}_j$, we may define a new covering of the attractor by means of the sets $T_\iii \mathbf{B}_j$. Since every $T_\iii$ is a similarity transformation each set $T_\iii \mathbf{B}_j$ is also a Euclidan ball, and its radius is smaller than that of $\mathbf{B}_j$ by a factor of exactly $r_\iii$. In this manner information about covers of $X$ by balls may be directly propagated to arbitrarily small scales. Indeed, by covering a self-similar set $X$ by a single ball $\mathbf{B}$ and then considering covers of the form $\{T_\iii \mathbf{B} \colon |\iii|=n\}$ it is an easy exercise to recover the upper bound $\dimh X\leq s$ directly from the definitions, where $s$ solves $\sum_{i=1}^N r_i^s=1$. This covering strategy is illustrated in Figure \ref{fi:covers} in the case of the Sierpi\'nski triangle. When the transformations $T_i$ are merely affine, however, Euclidean balls are no longer mapped to other Euclidean balls and the relationships between covers at different scales become more complicated: see the later Figure \ref{fi:affine-cover}.
 \subsection{Falconer's theorem on self-affine sets}
  In order to describe the central results obtained so far in the self-affine case, we will require some additional ideas and terminology. The following result from linear algebra will be fundamental in understanding the behaviour of covers of $X$ under the application of non-conformal affine transformations $T_\iii$, and hence in bounding the dimension of self-affine sets:
\begin{proposition}[Singular value decomposition]\label{pr:svd1}
Let $A \in M_d(\R)$. Then $A^TA$ has $d$ non-negative real eigenvalues $\sigma_1(A)^2\geq \cdots \geq \sigma_d(A)^2$, say, and there exists an orthonormal basis for $\R^d$ consisting of eigenvectors of $A^TA$. If $D \in M_d(\R)$ denotes the diagonal matrix with successive diagonal entries $\sigma_1(A),\ldots,\sigma_d(A)$ then there exist $U,V \in \O(d)$ such that $A=UDV^T$.
\end{proposition}
Here $\O(d)$ denotes the group of $d \times d$ matrices $U$ which preserve the Euclidean distance, or equivalently which satisfy $U^TU=I$. While the matrix $D$ in the decomposition $A=UDV^T$ is unique by definition, the matrices $U$ and $V$ are not: however, if $\hat{U}D\hat{V}^T$ is an alternative decomposition of $A$ then both $U^{-1}\hat{U}$ and $V^{-1}\hat{V}$ must preserve the eigenspaces of $D$ (see \cite[Theorem 2.6.5]{HoJo13}). The eigenvectors of $A^TA$ are in general called the \emph{singular vectors} of $A$.
We gather some useful properties of singular values below:
\begin{proposition}[Singular value identities and inequalities]\label{pr:svd2}
Let $A,B \in M_d(\R)$ and $1 \leq k \leq d$. Then:
\begin{equation}\label{eq:svd-obv}\sigma_1(A)=\|A\|,\qquad |\det A|=\prod_{i=1}^d \sigma_i(A),\qquad \sigma_k(A)=\sigma_k\left(A^T\right),\end{equation}
\begin{equation}\label{eq:svd-sub}\prod_{i=1}^k \sigma_i(AB) \leq \left(\prod_{i=1}^k \sigma_i(A)\right)\left(\prod_{i=1}^k \sigma_i(B)\right), \end{equation}
\begin{equation}\label{eq:svd-bounds}\sigma_k(A)\sigma_d(B) \leq \sigma_{k}(AB) \leq \sigma_1(A)\sigma_k(B),\end{equation}
and if $A$ is invertible then
\begin{equation}\label{eq:svd-dual}\sigma_k(A)=1/\sigma_{d+1-k}(A^{-1}).\end{equation}
\end{proposition}
A proof of Proposition \ref{pr:svd1} may be found in any textbook on matrix analysis worthy of the name. Most of the inequalities and identities presented in Proposition \ref{pr:svd2} are likewise standard: \eqref{eq:svd-obv} and \eqref{eq:svd-dual} are simple corollaries of Proposition \ref{pr:svd1} itself, and the bound $ \sigma_{k}(AB) \leq \sigma_1(A)\sigma_k(B)$ follows directly from any of several alternative characterisations of the singular values which we do not present here. The inequality $\sigma_k(A)\sigma_d(B) \leq \sigma_{k}(AB)$ is trivial when $B$ is non-invertible, and otherwise follows from the case $ \sigma_{d+1-k}(ABB^{-1}) \leq \sigma_1(AB)\sigma_{d+1-k}(B^{-1})$ of the second inequality in \eqref{eq:svd-bounds}, together with \eqref{eq:svd-dual}. The inequality \eqref{eq:svd-sub} is less commonly stated, but may be found for example as \cite[Theorem 3.3.4]{HoJo94}.

\begin{figure}
    \centering
    \subfloat[A self-affine set covered by a single ball $\mathbf{B}$.]{{\includegraphics[width=5.7cm]{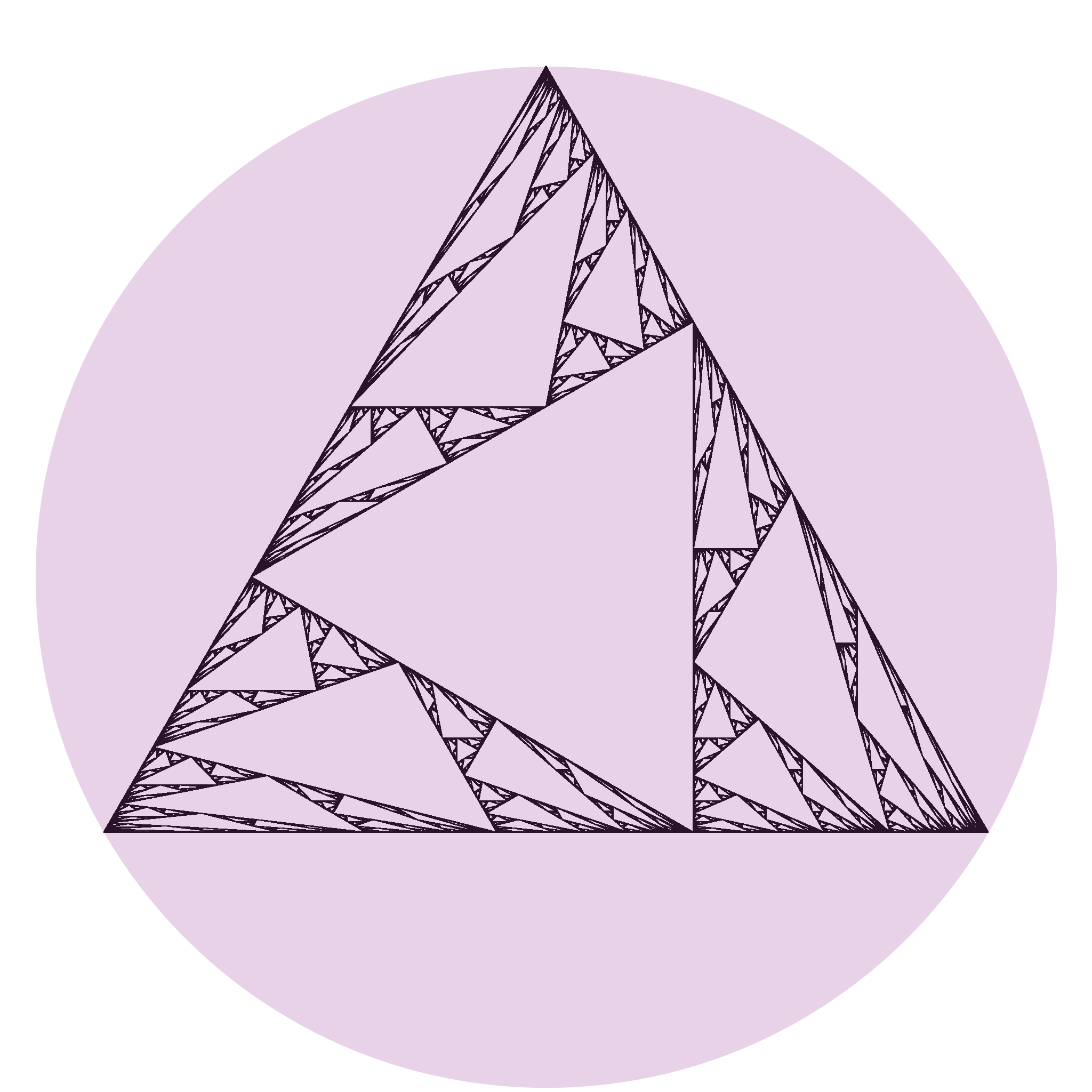}}}
    \quad
    \subfloat[The same  set covered by the images $T_\iii \mathbf{B}$ where $|\iii|=2$.]{{\includegraphics[width=5.7cm]{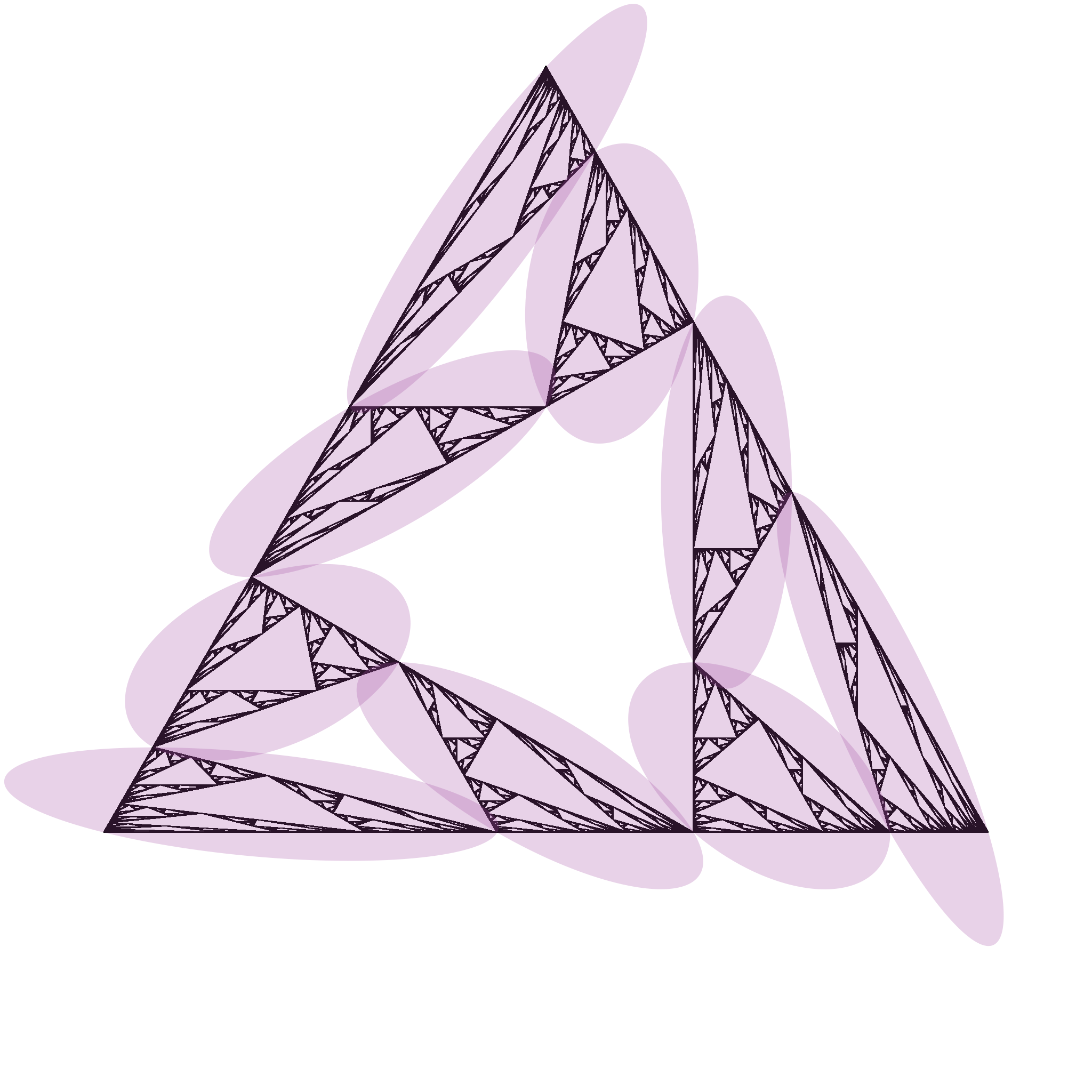}}}
    \caption{For self-affine sets which are not self-similar, na\"{\i}ve application of the covering strategy described \S\ref{ss:onepointone} and illustrated in Figure \ref{fi:covers} in general results in an inefficient cover using long, thin ellipsoids. To derive an efficient cover one dissects these ellipsoids into a larger number of smaller but rounder shapes.}
    \label{fi:affine-cover}%
\end{figure}

An immediate consequence of Proposition \ref{pr:svd1} is that if $\mathbf{B}$ is Euclidean ball in $\R^d$ of unit radius, then $A\mathbf{B}$ is an ellipsoid with semi-axes of length $\sigma_1(A),\ldots,\sigma_d(A)$. It is precisely this fact which underlies the standard upper bound for the dimension of a self-affine set $X=\bigcup_{i=1}^N T_iX$. Beginning with a  Euclidean ball $\mathbf{B}$ which contains the attractor, one considers a covering of $X$ by sets of the form $T_\iii \mathbf{B}$. Each such set is an ellipsoid whose axes have lengths respectively $\sigma_1(A_\iii),\ldots,\sigma_d(A_\iii)$ times the diameter of $\mathbf{B}$. Now, the definition of the Hausdorff dimension perceives only the diameter of each set in a given cover, regardless of its shape; one should not, therefore, expect covers by elongated shapes to be efficient, since a cover by thin ellipsoids of length $\delta$ (say) can in general be replaced with a cover using a much smaller number of balls of the same diameter $\delta$. This problem is illustrated in Figure \ref{fi:affine-cover}. 

Given an elongated covering set, then, one is presented with various options as to how to efficiently cover it with rounder shapes such as spheres or cubes. A cuboid of dimensions $1\times 4\times 9$, for example, may be covered with a single cube of side length $9$, or with three cubes of side length $4$, or with $36$ cubes of side length $1$. In general, which option yields the best bound for the quantity $\mathcal{H}^s_\delta$  depends on the dimension $s$ being considered. In practice, to cover an ellipsoid having semi-axes $\sigma_1(A_\iii),\ldots,\sigma_d(A_\iii)$ with sets $U_j$ in a way which makes $\sum_j (\diam U_j)^s$ small it turns out to be most efficient to use a covering by cubes whose diameter matches the parameter $\sigma_{\lceil s\rceil}(A_\iii)$. The number of cubes needed is then approximately $(\sigma_1(A_\iii)\cdots \sigma_{\lceil s\rceil}(A_\iii)) / \sigma_{\lceil s\rceil}(A_\iii)^{\lceil s\rceil}$, and the total contribution $\sum_j(\diam U_j)^s$ of the diameters of the cubes used to cover this ellipsoid is proportional to the quantity
\[\left(\frac{\prod_{j=1}^{\lceil s\rceil} \sigma_j(A_\iii)}{\sigma_{\lceil s\rceil}(A_\iii)^{\lceil s\rceil}}\right)\sigma_{\lceil s\rceil}(A_\iii)^s = \sigma_1(A_\iii)\cdots \sigma_{\lfloor s\rfloor}(A_\iii) \sigma_{\lceil s\rceil}(A_\iii)^{s-\lfloor s\rfloor}.\]
Given $A \in M_d(\R)$, then, for each $s \in [0,d]$ define
\[\varphi^s(A):=\sigma_1(A)\sigma_2(A)\cdots \sigma_{\lfloor s\rfloor}(A) \sigma_{\lceil s\rceil}(A)^{s-\lfloor s\rfloor}.\]
It is an easy consequence of \eqref{eq:svd-bounds} that $\varphi^s(AB)\leq \varphi^s(A)\varphi^s(B)$ for all $A,B \in M_d(\R)$. Given $\A=(A_1,\ldots,A_N) \in \GL_d(\R)$, let us say that $\A$ is \emph{contracting} if there exists a norm $\threebar{\cdot}$ on $\R^d$ such that $\max_i \threebar{A_i}<1$, and \emph{strongly contracting} if there exists a norm such that $\max_{1 \leq i<j\leq N}\threebar{A_i}+\threebar{A_j}<1$. Given a contracting $\A\in \GL_d(\R)^N$ and a tuple $\mathbf{v}=(v_1,\ldots,v_N)\in (\R^d)^N$ we define an iterated function system $(T_1^{\mathbf{v}},\ldots,T_N^{\mathbf{v}})$ by $T_i^{\mathbf{v}}:=A_ix+v_i$ for every $i=1,\ldots,N$ and $x \in \R^d$. In this case we say that $\A$ is the \emph{linearisation} of the iterated function system $(T_1^{\mathbf{v}},\ldots,T_N^{\mathbf{v}})$. The general theory of self-affine fractals was initiated by the following landmark result of K. Falconer in the 1988 article \cite{Fa88}: 
\begin{theorem}[Falconer]\label{th:falconer}
Let $\A =(A_1,\ldots,A_N) \in \GL_d(\R)^N$. Then:
\begin{enumerate}[(i)]
\item
If $\A$ is contracting, then for every $\mathbf{v} \in (\R^d)^N$ the attractor $X$ of the iterated function system $(T_1^{\mathbf{v}},\ldots,T_N^{\mathbf{v}})$  satisfies
\[\dimh X \leq \dimaff \A:=\inf\left\{s \in [0,d]\colon \sum_{n=1}^\infty \sum_{|\iii|=n} \varphi^s(A_\iii)<\infty\right\}.\]
\item
If $\A$ is strongly contracting, then for Lebesgue almost every $\mathbf{v} \in (\R^d)^N$ the attractor $X$ of the iterated function system $(T_1^{\mathbf{v}},\ldots,T_N^{\mathbf{v}})$  satisfies $\dimh X = \dimaff \A$.
\end{enumerate}
Here $\dimaff \A$ is called the \emph{affinity dimension} of $\A$, and is understood to equal $d$ in the case where $ \sum_{n=1}^\infty \sum_{|\iii|=n} \varphi^s(A_\iii)$ diverges for every $s \in [0,d]$.
\end{theorem}
Here we present Falconer's hypotheses in a slightly updated form: in the original work \cite{Fa88} the contraction hypothesis was defined only with reference to the Euclidean norm, and the strong contraction hypothesis was the more restrictive condition $\max_i \|A_i\|<\frac{1}{3}$. The version stated above corresponds more closely with that presented in the recent textbook \cite{BaSiSo23}. The strong contraction hypothesis cannot be substantially relaxed: some examples are surveyed in \cite{Mo23}.

If every $T_i$ is a similarity transformation with contraction ratio $r_i$ then the affinity dimension coincides with the solution to $\sum_{i=1}^N r_i^s=1$ as in Hutchinson's theorem. In contrast to Hutchinson's theorem, however, Falconer's expression for the dimension of the attractor is shown to hold for Lebesgue almost every translation component $\mathbf{v}$ as opposed to in the case of an explicit condition on the separation of the images $T_iX$. The problem of realising Falconer's dimension formula under an explicit condition  analogous to that in Hutchinson's theorem has been the subject of extensive research. Surprisingly, it is easier to construct non-examples than examples: for instance, the Cartesian product of two self-affine sets is always self-affine, but very rarely has Hausdorff dimension equal to $\dimaff \A$. It is an interesting exercise for the reader to demonstrate that the product of two middle-interval Cantor sets has dimension smaller than the affinity dimension provided by its natural representation as a self-affine set. The following open question is thus the central problem of self-affine fractal geometry:
\begin{question}[Dimension problem for self-affine sets]
Let $(T_1,\ldots,T_N)$ be an affine iterated function system with linearisation $\A$ and attractor $X$. Under what hypotheses is it the case that $\dimh X= \dimaff \A$?
\end{question}
Let us say that $(T_1,\ldots,T_N)$ satisfies the \emph{open set condition} if there exists a nonempty open set $U \subset \R^d$ such that the sets $T_iU$ are pairwise disjoint subsets of $U$, and let us say that the \emph{strong open set condition} holds if  the set $U$ can be  chosen so that additionally $U \cap X \neq \emptyset$. We will also say that $\A \in \GL_d(\R)^N$ is \emph{irreducible} if there is no proper nonzero subspace of $\R^d$ which is preserved by every $A_i$, and \emph{strongly irreducible} if no finite union of proper nonzero subspaces is preserved by every $A_i$. A tuple which is not irreducible will be called \emph{reducible}. It is known by work of Bedford and McMullen \cite{Be84,Mc84} that the strong open set condition does not guarantee the conclusion $\dimh X=\dimaff \A$ in the reducible case, nor (by work of Fraser in \cite{Fr12}) does this implication hold in the irreducible but not strongly irreducible case. By an example of Edgar \cite{Ed92} the open set condition does not imply $\dimh X=\dimaff \A$ even in the strongly irreducible case, since the open set condition can hold even in cases where the attractor is a singleton. A more recent programme of research -- see \cite{BaHoRa19,MoSe23,Ra24} -- has demonstrated that in dimension $d \leq 3$, strong irreducibility and the strong open set condition together imply the conclusion $\dimh X=\dimaff \A$. In dimension $d \geq 4$ the problem remains wide open.

\subsection{Projections of fractal sets and an extension of Falconer's theorem}

Let us now turn our attention to the dimensions of \emph{projections} of fractal sets. If $X \subset \R^d$ is arbitrary, and if $Q \in M_d(\R)$ is an orthogonal projection, it is an elementary exercise to show that $\dimh QX \leq \min \{\dimh X, \rank Q\}$. Let us say that $Q$ is a \emph{full} projection of $X$ if $\dimh QX$ is equal to this upper bound, and an \emph{exceptional} projection otherwise. It is easily seen that if $X$ is an open subset of a $k$-dimensional linear subspace $V$, say, then projections whose kernel intersects $V$ are exceptional and all other projections are full. A classical result of J.M Marstrand \cite{Ma54}, substantially extended and generalised by P. Mattila in \cite{Ma75}, extends this picture to arbitrary Borel sets $X$:
\begin{theorem}[Marstrand's theorem]
Let $X \subset \R^d$ be Borel, let $1 \leq k\leq d$ and let $0<s<\dimh X$. Then the set of all orthogonal projections $Q \in M_d(\R)$ which have rank precisely $k$ and which satisfy $\dimh QX \leq s$ has Hausdorff dimension at most $k(d-k)-(\dimh X-s)$.
\end{theorem}
The set of all rank-$k$ orthogonal projections is a smooth submanifold of $M_d(\R)$ of dimension $k(d-k)$, so this result guarantees in particular that the set of exceptional projections of a given set $X$ always has Lebesgue measure zero. On the other hand the precise structure of this measure-zero set remains rather unclear. The preceding example --  in which $X$ is an open subset of a lower-dimensional vector space -- demonstrates that the exceptional set can contain every projection whose kernel intersects a given linear subspace. A remarkable programme of research initiated by F\"assler and Orponen -- see \cite{FaOr14,GaGuGuHaMaWa22,KaOrVe25,PrYaZa22} -- shows that in the case $X\subset \R^3$, this is the \emph{only} kind of $C^2$ submanifold which can be contained in the exceptional set. In higher dimensions it turns out that broader classes of exceptional projections are possible, so that for example the set of exceptional rank-$2$ projections of a set $X\subset \R^4$ can contain smooth manifolds which do not conform to the preceding description. In this article we will describe and present some examples of this phenomenon in the particular case where $X$ is a self-affine set.

At the core of the examples which we present is the following recent theorem of D.-J. Feng, Yu-Hao Xie, \cagri Sert and the author, which strictly extends Theorem \ref{th:falconer}. We retain the notational conventions of that theorem:
\begin{theorem}[Feng-Xie \cite{FeXi25}, Morris-Sert \cite{MoSe25}]\label{th:fxms}
Let $\A=(A_1,\ldots,A_N) \in \GL_d(\R)^N$ and $Q \in M_d(\R)$. Then:
\begin{enumerate}[(i)]
\item
If $\A$ is contracting, then for every $\mathbf{v} \in (\R^d)^N$ the attractor $X$ of the iterated function system $(T_1^{\mathbf{v}},\ldots,T_N^{\mathbf{v}})$  satisfies
\[\dimh QX \leq \dimaffQ\A:=\inf\left\{s \in [0,\rank Q] \colon \sum_{n=1}^\infty \sum_{|\iii|=n} \varphi^s(QA_\iii)<\infty\right\}.\]
\item
If $\A$ is strongly contracting, then for Lebesgue almost every $\mathbf{v} \in (\R^d)^N$ the attractor $X$ of the iterated function system $(T_1^{\mathbf{v}},\ldots,T_N^{\mathbf{v}})$  satisfies $\dimh QX = \dimaffQ \A$.
\end{enumerate}
Here we similarly define $\dimaffQ \A:=\rank Q$ if the series diverges for every $s \in [0,\rank Q]$.
\end{theorem}
The significance of this result is that the inequality $\dimaffQ \A<\dimaff \A$ is possible across a broad range of examples, which are substantially characterised in \cite{MoSe25}. It is further shown in that work that the set $\{Q \in M_d(\R)\colon \dimaffQ \A<\dimaff \A\}$ always has the structure of an algebraic variety which is invariant under right-multiplication by the matrices $A_i$. In such cases one may therefore exhibit self-affine sets $X$ which admit nontrivial manifolds of exceptional projections.

The proof of Theorem \ref{th:fxms} is long, and we will not attempt to review the whole of it here. However, the proof of the \emph{upper} bound in particular is remarkably tractable. In the remainder of this article we provide an exposition of the proof of the upper bound in Theorem \ref{th:fxms}, and we will go on to present some examples in which $\dimaffQ \A<\dimaff \A$ for an explicitly-identifiable family of projections $Q$. This in particular exhibits some manifolds which can occur as subsets of the exceptional set in Marstrand's Theorem. 

We close this introduction with some very brief remarks on the lower bound. In all of the theorems stated in this introduction, lower bounds may be established via the \emph{energy method}, as follows. Suppose that one is able to construct a measure $\mu$ supported on $X$ such that  $\iint \|x-y\|^{-s}d\mu(x)d\mu(y)<\infty$. Let $\mathbf{B}(x,r)$ denote the closed Euclidean ball of centre $x$ and radius $r$. If the integral converges then there exist a positive-measure set $E\subset X$ and $C>0$ such that $\int  \|x-y\|^{-s}d\mu(y) \leq C$ for all $x \in E$, hence by Markov's inequality $\mu(\mathbf{B}(x,r)) \leq Cr^s$ for all $x \in E$ and $r>0$. If $(U_j)_{j=1}^\infty$ is any cover of $E$ by sets of diameter at most $\delta>0$, take $(x_j)_{j=1}^\infty$ to be a sequence of points in $E$ such that $x_j \in U_j$ for all $j \geq 1$ and note that $U_j \subseteq \mathbf{B}(x_j, \diam U_j)$ for every $j$. For every $t \in (0,s)$ we have
\begin{align*}0<\mu(E)\leq \sum_{j=1}^\infty \mu(U_j) &\leq \sum_{j=1}^\infty \mathbf{B}(x_j, \diam U_j)\\
&\leq C\sum_{j=1}^\infty (\diam U_j)^s \leq C\delta^{s-t}\sum_{j=1}^\infty (\diam U_j)^t\end{align*}
and therefore
\[\mathcal{H}^t_\delta (X) \geq \mathcal{H}^t_\delta(E) \geq \frac{\mu(E)}{C\delta^{s-t}}\]
so that $\lim_{\delta \to 0} \mathcal{H}^t_\delta(X)=+\infty$ and consequently $\dimh X \geq t$. The lower bound in Hutchinson's theorem constructs such a measure $\mu$ in a rather direct manner via the characterisation $\mu(T_\iii X)=r_\iii$ for every word $\iii$. Falconer's theorem instead proceeds by defining a measure $\mu_{\mathbf{v}}$ which depends on the parameter $\mathbf{v} \in (\R^d)^N$, and evaluating a triple integral over $x$, $y$ and $\mathbf{v}$ to show that for Lebesgue a.e. $\mathbf{v}$ the relevant double integral in $x,y$ converges. The proof of Theorem \ref{th:fxms} is made more difficult by the fact that the most natural choice of $\mu$ in general results in a divergent integral. To circumvent this problem one in effect chooses instead a set $E\subset X$ on which the double integral converges, although the proofs in \cite{FeXi25,MoSe25} are not expressed in quite those terms. 

\section{Proof of the upper dimension bound}
We begin our exposition by presenting a proof of the upper bound in Theorem \ref{th:fxms}. The proof of this bound is a surprisingly straightforward adaptation of Falconer's original proof of Theorem \ref{th:falconer}:
\begin{proposition}\label{pr:ubend}
Let $X\subset \R^d$ be the attractor of an affine iterated function system with linearisation $\A \in \GL_d(\R)^N$, let $Q \in M_d(\R)$, and define
\[\dimaffQ \A:= \inf\left\{s \in [0,\rank Q]\colon\sum_{n=1}^\infty \sum_{|\iii|=n} \varphi^s(QA_\iii)<\infty\right\} \]
where this quantity is understood to equal $\rank Q$ if the set is empty. Then $\dimh QX \leq \dimaffQ \A$.
\end{proposition}
At some cost to brevity this result may be easily extended so as to give the same bound for the upper box dimension of $QX$, but we do not pursue this here: see \cite{MoSe25}.
\begin{proof}
We adapt the arguments of Falconer's upper bound from the classic article \cite{Fa88}. If $\dimaffQ \A=\rank Q$ then the desired bound is trivial. Otherwise, fix an arbitrary closed Euclidean ball $\mathbf{B}$ which contains $X$, choose an arbitrary $s \in [0,\rank Q]$ such that the series $\sum_{n=1}^\infty\sum_{|\iii|=n}\varphi^s(QA_\iii)$ converges, and let $\delta>0$ also be arbitrary. By equivalence of norms there exists $C\geq 1$ such that $C^{-1}\|B\|\leq \threebar{B}\leq C\|B\|$ for every $B \in M_d(\R)$, and since $\A$ is contracting it follows that we may choose an integer $m(\delta) \geq 1$ large enough that $\max_{|\iii|=m(\delta)} \|A_\iii\|<\delta/\sqrt{d}$. 
By iterating the equation $\bigcup_{i=1}^N T_iX$ we obtain $X=\bigcup_{|\iii|=m(\delta)} T_\iii X$, and therefore
\[QX =\bigcup_{|\iii|=m} QT_\iii X \subseteq \bigcup_{|\iii|=m(\delta)} QT_\iii \mathbf{B}.\] 

Consider any one of the sets $QT_\iii \mathbf{B}$ where $|\iii|=m(\delta)$. Clearly this set is a translated copy of the set $QA_\iii \mathbf{B}$, and the latter is a (possibly degenerate) closed ellipsoid whose axes have length precisely $\sigma_j(QA_\iii)\cdot \diam \mathbf{B}$ for $j=1,\ldots,d$. In particular $T_\iii Q\mathbf{B}$ is contained in a (possibly degenerate) closed cuboid whose side lengths are $\sigma_j(QA_\iii)\cdot \diam \mathbf{B}$ for $j=1,\ldots,d$. Let us cover the latter with closed cubes of side length $\sigma_{\lceil s\rceil}(QA_\iii)<\delta/\sqrt{d}$. (This number is nonzero since $QA_\iii$ has rank $k \geq \lceil s\rceil$.) Dividing each axis of the cuboid into equal intervals of length $\sigma_{\lceil s\rceil}(QA_\iii)$ together with at most one ``remainder'' interval of length less than $\sigma_{\lceil s\rceil}(QA_\iii)$, and using the resulting cubic mesh to cover the cuboid,  it is not difficult to see that the number of cubes required to cover the cuboid is not greater than
\begin{align*}\prod_{\ell=1}^d \left(\left\lfloor \frac{\sigma_\ell (QA_{\iii}) \diam \mathbf{B}}{\sigma_{\lceil s\rceil}(QA_{\iii})} \right\rfloor +1 \right) &
\leq\left(2\left\lceil\diam \mathbf{B}\right\rceil\right)^d \cdot  \prod_{\ell=1}^d \max\left\{1, \frac{\sigma_\ell (QA_{\iii})}{\sigma_{\lceil s\rceil}(QA_{\iii})}\right\} \\
&= \left(2\left\lceil\diam \mathbf{B}\right\rceil\right)^d  \prod_{\ell=1}^{\lceil s\rceil} \frac{\sigma_\ell (QA_{\iii})}{\sigma_{\lceil s\rceil}(QA_{\iii})} \\
& =  \left(2\left\lceil\diam \mathbf{B}\right\rceil\right)^d  \cdot \frac{ \varphi^s (QA_{\iii})}{\sigma_{\lceil s\rceil}(QA_{\iii})^s}.\end{align*}
(Here we have used the inequality $\lfloor xy\rfloor+1 \leq 2\lceil x\rceil \max\{1,y\}$  which holds for all $x>0$ and $y\geq 0$ and which is easily verified by the reader.) Thus for every $\iii$ of length $m(\delta)$ we may choose finitely many cubes $C_1(\iii),\ldots,C_{k(\iii,\delta)}(\iii)$ of diameter $\sqrt{d}\cdot \sigma_{\lceil s\rceil}(QA_\iii)<\delta$ whose union covers $T_\iii \mathbf{B}$ and such that $\sum_{j=1}^{k(\iii,\delta)} (\diam C_j(\iii))^s \leq (2\sqrt{d}\lceil \diam \mathbf{B}\rceil)^d \varphi^s(QA_\iii)$. Applying this strategy to every word $\iii$ of length $m(\delta)$ yields the bound
\begin{align*}\mathcal{H}^s_\delta(QX) &\leq \mathcal{H}^s_\delta\left(\bigcup_{|\iii|=m(\delta)} QT_\iii \mathbf{B}\right)\\
& \leq \sum_{|\iii|=m(\delta)} \sum_{j=1}^{k(\iii,\delta)} (\diam C_j(\iii))^s \\
& \leq \left(2\sqrt{d}\lceil \diam \mathbf{B}\rceil\right)^d\cdot \sum_{|\iii|=m(\delta)} \varphi^s(QA_\iii)\\
& \leq \left(2\sqrt{d}\lceil \diam \mathbf{B}\rceil\right)^d\left(\sum_{n=1}^\infty \sum_{|\iii|=n} \varphi^s(QA_\iii)\right)<\infty\end{align*}
where we note that the last quantity is independent of the choice of $\delta>0$. The bound $\dimh QX \leq s$ follows.\end{proof}
The quantity $\dimaffQ \A$ as defined above is \emph{a priori} rather difficult to estimate, so to simplify its analysis we make a further definition:
\begin{definition}
Let $\A \in \GL_d(\R)^N$, $Q \in M_d(\R)$ and $s \in [0,\rank Q]$. We define the \emph{pressure of $\A$ at $s$} to be the real number
\[P(\A,s):=\lim_{n \to \infty} \frac{1}{n}\log \sum_{|\iii|=n}\varphi^s(A_\iii),\]
and the \emph{$Q$-projected pressure of $\A$ at $s$} to be the quantity
\[P_Q(\A,s):=\limsup_{n \to \infty} \frac{1}{n}\log \sum_{|\iii|=n} \varphi^s(QA_\iii).\]
\end{definition}
The quantity $P(\A,s)$ was considered implicitly by Falconer in \cite[\S4]{Fa88} and in a more explicit form by K\"aenm\"aki in \cite{Ka04}. The existence of this limit is guaranteed by subaddivity, as a consequence of the general inequality $\varphi^s(AB)\leq \varphi^s(A)\varphi^s(B)$. 

The inequality $P_Q(\A,s) \leq P(\A,s)$ is obvious, so in particular $P_Q(\A,s)$ is never equal to $+\infty$. Less obviously, this pressure also never equals $-\infty$. To see this, suppose that $\rank Q=k \geq s$. Since $\sigma_j(QB) \geq \sigma_k(QB) \geq \sigma_k(Q)\sigma_d(B)$ for all $j\leq k$ and $B \in \GL_d(\R)$,
\begin{align*}\limsup_{n \to \infty} \frac{1}{n}\log \sum_{|\iii|=n} \varphi^s(QA_\iii) &\geq \liminf_{n \to \infty}\frac{1}{n} \log \sum_{|\iii|=n} \sigma_k(Q)^s \sigma_d(A_\iii)^s\\
& =\lim_{n \to \infty} \frac{1}{n}\log \sum_{|\iii|=n} \sigma_d(A_\iii)^s\geq \log \sum_{i=1}^N \sigma_d(A_i)^s>-\infty.\end{align*}
Here we have used the fact that the sequence $n \mapsto -\log \sum_{|\iii|=n}\sigma_d(A_\iii)^s$ is subadditive as a consequence of the general inequality $\sigma_d(AB) \geq \sigma_d(A)\sigma_d(B)$, which itself follows from Proposition \ref{pr:svd2}. The limit superior which defines $P_Q(\A,s)$ is in fact a limit (see \cite{FeXi25,MoSe25}) but the proof of this result is rather complicated and will not be necessary for our exposition.

The utility of the above definition in understanding $\dimaffQ \A$ lies in the following result, which also has roots in Falconer's article \cite{Fa88}:
\begin{proposition}\label{pr:cty}
Let $\A \in \GL_d(\R)^N$ be contracting and let $Q \in M_d(\R)$. Then the function 
$s\mapsto P_Q(\A,s)$ is continuous and strictly decreasing on $[0,\rank Q]$ and is convex on every interval $[k-1,k]\subseteq [0,\rank Q]$ with $k$ integer. Additionally
\[\dimaffQ \A=\inf\left\{s \in [0,\rank Q] \colon P_Q(\A,s) \leq 0\right\}\]
when the latter set is nonempty, and $\dimaffQ \A=\rank Q$ otherwise.
\end{proposition}
\begin{proof}
Define $e^{-\alpha}:=\max_{1 \leq i \leq N} \threebar{A_i}$ and $e^{-\beta}:=\min_{1 \leq i \leq N} \sigma_d(A_i)$, both of which lie in $(0,1)$. By equivalence of norms we may choose $C \geq 1$ such that $C^{-1}\|B\| \leq \threebar{B}\leq C\|B\|$ for every $B \in M_d(\R)$.  If $0 \leq k-1 \leq s < s+t \leq  k \leq \rank Q$, then for every $n \geq 1$
\begin{align*}\sum_{|\iii|=n}\varphi^{s+t}(QA_\iii) &\leq\left(\max_{|\iii|=n} \sigma_k(QA_\iii)^t\right) \sum_{|\iii|=n}\varphi^{s}(QA_\iii)  \\
&\leq\left(C^t\threebar{Q}^t\cdot \max_{|\iii|=n} \threebar{A_\iii}^t\right) \sum_{|\iii|=n}\varphi^{s}(QA_\iii)\\
& \leq\left(C^t\threebar{Q}^t e^{-n\alpha t}\right) \sum_{|\iii|=n}\varphi^{s}(QA_\iii)  \end{align*}
and 
\begin{align*}\sum_{|\iii|=n}\varphi^{s+t}(QA_\iii)& \geq\left(\min_{|\iii|=n} \sigma_k(QA_\iii)^t\right) \sum_{|\iii|=n}\varphi^{s}(QA_\iii)  \\
&\geq\left(\sigma_k(Q)^t \cdot \min_{|\iii|=n} \sigma_d(A_\iii)^t\right) \sum_{|\iii|=n}\varphi^{s}(QA_\iii)\\
&\geq\sigma_k(Q)^te^{-n\beta t} \sum_{|\iii|=n}\varphi^{s}(QA_\iii)\end{align*}
using Proposition \ref{pr:svd2}. It follows directly that
\[P_Q(\A,s)-\beta t\leq P_Q(\A,s+t) \leq P_Q(\A,s) -\alpha t,\]
and this proves that $s\mapsto P_Q(\A,s)$ is Lipschitz continuous and strictly decreasing on the interval $[k-1,k]$. To prove convexity on this interval we simply note that for every $\alpha \in (0,1)$ and $n \geq 1$
\begin{align*}\sum_{|\iii|=n}\varphi^{(1-\alpha)s + \alpha(s+t)}(QA_\iii) &=\sum_{|\iii|=n} \varphi^s(QA_\iii)^{1-\alpha}  \varphi^{s+t}(QA_\iii)^\alpha\\
&\leq\left(\sum_{|\iii|=n} \varphi^s(QA_\iii)\right)^{1-\alpha}  \left(\sum_{|\iii|=n} \varphi^s(QA_\iii)\right)^{\alpha}   \end{align*}
by H\"older's inequality with $p=\frac{1}{1-\alpha}$, $q=\frac{1}{\alpha}$, and the result follows.
Since $k$ was arbitrary this completes the proof.
\end{proof}
Combining Propositions \ref{pr:ubend} and \ref{pr:cty} we arrive at the following tool for detecting exceptional projections of self-affine sets:
\begin{corollary}\label{co:cond}
Let $X\subset \R^d$ be the attractor of an iterated function system with linearisation $\A \in \GL_d(\R)^N$, let $Q \in M_d(\R)$ and suppose that $\dimaff \A \in [0,\rank Q]$. If $P_Q(\A,\dimaff \A)<0$, then $\dimh QX<\dimaff \A$.
\end{corollary}
\begin{proof}
Since clearly $P_Q(\A,0) =\log N>0$, it follows easily from Proposition \ref{pr:cty} that $s \mapsto P_Q(\A,s)$ has a unique zero in $(0,\dimaff \A)$. The series $\sum_{n=1}^\infty \sum_{|\iii|=n}\varphi^s(QA_\iii)$ obviously converges when $P_Q(\A,s)$ is strictly negative and diverges when it is strictly positive, so the unique zero of $P_Q(\A,s)$ must be precisely $\dimaffQ \A$. Thus $\dimaffQ \A<\dimaff \A$. The dimension bound follows by Proposition \ref{pr:ubend}. \end{proof}

\section{Construction of examples}
\subsection{Underlying principles and preliminaries}

Corollary \ref{co:cond} above gives us a crisp criterion with which to construct examples of self-affine sets with exceptional projections: to do this, we must find examples of $\A$ and $Q$ such that $P_Q(\A,s)<P(\A,s)$ for $s=\dimaff \A$ and for some nonzero $Q \in M_d(\R)$. On the other hand, it is not immediately clear how such examples might be constructed.

The following simple example serves as a guide. Let $(B_1,\ldots,B_N) \in \GL_{d_1}(\R)^N$ and $(B_1',\ldots,B_N') \in \GL_{d_2}(\R)^N$, and consider the tuple $\A=(A_1,\ldots,A_N)\in \GL_{d_1+d_2}(\R)^N$ defined by
\begin{equation}\label{eq:A}A_i=\begin{pmatrix} B_i &0\\0&B_i'\end{pmatrix}.\end{equation}
The block diagonal form of each matrix $A_\iii$ implies that the list of singular values $\sigma_1(A_\iii),\ldots,\sigma_{d_1+d_2}(A_\iii)$ consists simply of the singular values of $B_\iii$ and of $B_\iii'$, intermingled with one another in some \emph{a priori} unknown order. To simplify this situation we impose an additional condition. Suppose that for some $\kappa>0$
\begin{equation}\label{eq:weak-dom}\max_{1 \leq i\leq N} \sigma_1(B_i') < e^{-\kappa}\min_{1 \leq i \leq N} \sigma_d(B_i).\end{equation}
Since $\sigma_1$ is submultiplicative and $\sigma_d$ supermultiplicative we find that
\[\sigma_{1}(B_\iii') < e^{-\kappa |\iii|} \sigma_{d_1}(B_\iii)\]
for every word $\iii$; consequently the largest $d_1$ singular values of $A_\iii$ must be precisely the $d_1$ singular values of $B_\iii$, and the smallest $d_2$ singular values of $A_\iii$  the singular values of $B_\iii'$. Under the hypothesis \eqref{eq:weak-dom} we thus have
\[\sigma_{d_1+1}(A_\iii) \leq \sigma_1(B_\iii') <e^{-\kappa |\iii|}\sigma_{d_1}(B_\iii) \leq e^{-\kappa |\iii|}\sigma_{d_1}(A_\iii)\]
for every finite word $\iii$. It follows that if $Q \in M_d(\R)$ is given by 
\[Q=\begin{pmatrix} 0&0\\ 0&I\end{pmatrix}\]
where the dimensions of the block matrices are the same as in \eqref{eq:A}, then for every $k=1,\ldots,d_1$ and every word $\iii$ 
\[\sigma_k(QA_\iii)=\sigma_k(B_\iii')\leq \sigma_1(B_\iii') <e^{-\kappa |\iii|}\sigma_{d_1}(A_\iii) \leq e^{-\kappa |\iii|} \sigma_k(B_\iii)=\sigma_k(A_\iii).\] 
Hence $P_Q(\A,s)\leq P(\A,s)-s\kappa<P(\A,s)$ for all $s \in (0,d_1]$, allowing examples of self-affine sets for which $Q$ is an exceptional projection to be constructed.

The above example goes through with minor modifications in the more general case where the upper-right block of every $A_i$ is allowed to be an arbitrary $d_1 \times d_2$ matrix, in which case the key inequalities above are satisfied up to multiplication by an additional sub-exponential error term. What is less immediately obvious is the following further extension: if $Q$ again has rank $d_1$ but its  kernel intersects the span of the first $d_1$ co-ordinate basis vectors on a subspace of dimension $\ell$, say, then the pressure again drops after composition with $Q$. In this case, the first $(d_1-\ell)$ singular values of $QA_\iii$ arise from the upper-left block of $A_\iii$ and the remaining nonzero values arise from the lower-right block, with the result that for all $s \in (d_1-\ell, d_1]$ we again have $P_Q(\A,s)<P(\A,s)$. Since there are in general at least an $\ell(d_1-\ell)$-dimensional family of orthogonal projections $Q$ with this property, this shows in particular that the class of exceptional projections can be rather large.

These examples suggest the following more general strategy. Consider a tuple $\A$ for which there exists $\kappa>0$ and an integer $k$ such that
\[\frac{\sigma_{k+1}(A_\iii)}{\sigma_k(A_\iii)} <e^{-\kappa|\iii|}\]
for all long enough words $\iii$, a property which is called \emph{$k$-domination} (see \cite{BoGo09}). Suppose also that there exists a projection $Q$ of rank $k$ such that for every finite word $\iii$, the kernel of $Q$ intersects the span of the first $k$ singular vectors of $A_\iii$ in a subspace of dimension at least $\ell \geq 1$. Then the largest $k-\ell$ singular values of $QA_\iii$ are bounded above by the first $k-\ell$ singular values of $A_\iii$, whereas the next $\ell$ singular values of $QA_\iii$ will be exponentially smaller than the corresponding singular values of $A_\iii$; consequently $P_Q(\A,s)<P(\A,s)$ for all $s \in (k-\ell,k]$. In particular if $\dimaff \A \in (k-\ell,k]$ then a dimension drop occurs in the $Q$-projected image of the attractor of a typical iterated function system with linearisation $\A$. We will see that unlike the block diagonal construction sketched above, this more general situation is compatible with $\A$ being irreducible or even strongly irreducible.  In the next two subsections we will give constructions of examples which conform to this general schema.

Before presenting the first of our examples, we recall some general facts concerning tensor products of matrices. Given two square or rectangular matrices $A=[a_{ij}]$ and $B=[b_{ij}]$ with respective dimensions $n_1\times m_1$ and $n_2 \times m_2$, their Kronecker tensor product $A\otimes B$ is defined to be the $n_1n_2 \times m_1m_2$ matrix
\[A\otimes B:=\begin{pmatrix}
a_{11}B & a_{12}B &\cdots &a_{1m_1}B\\ 
a_{21}B & a_{22}B &\cdots &a_{2m_1}B\\
\vdots &\vdots & \ddots &\vdots\\ 
a_{n_1 1}B &a_{n_12}B & \cdots &a_{n_1m_1}B  
\end{pmatrix}.\]
The fundamental properties of this construction are summarised as follows:
\begin{proposition}
Let $A, A_1, A_2 \in M_{d_1}(\R)$,  $B, B_1, B_2 \in M_{d_2}(\R)$ and $\lambda \in \R$. Then:
\[(A_1 \otimes B_1)(A_2\otimes B_2) = (A_1A_2)\otimes (B_1B_2)\]
\[\lambda(A \otimes B) = (\lambda A)\otimes B = A\otimes (\lambda B)\]
and
\[(A\otimes B)^T = A^T \otimes B^T,\]
and the singular values of $A \otimes B$ are precisely the $d_1d_2$ numbers of the form $\sigma_i(A)\sigma_j(B)$, where $1 \leq i \leq d_1$ and $1 \leq j \leq d_2$. In particular,
\[\|A\otimes B\| = \|A\|\cdot \|B\|.\]
\end{proposition}
A proof of the above properties may be found in  \cite[\S4]{HoJo94}. It follows from the above result that if $Q$ is an orthogonal projection then so are $I\otimes Q$ and $Q\otimes I$. The operation $\otimes$ is associative: given $A \in M_d(\R)$ and $k \geq 1$, we let $A^{\otimes k} \in M_{d^k}(\R)$ denote the tensor product $A \otimes A \otimes \cdots \otimes A$ of $k$ copies of $A$, which in view of the above proposition has norm precisely $\|A\|^k$ and satisfies the identity $(AB)^{\otimes k}\equiv A^{\otimes k} B^{\otimes k}$ for $A,B \in M_{d}(\R)$.

\subsection{A concrete example constructed via Kronecker products}

We begin this tour by noting the following very specific example, which dates back to the very earliest origins of this programme of research. 
\begin{proposition}\label{pr:old}
Define $\A:=(A_1, A_2) \in \GL_4(\R)^2$ by
\[A_1:= \frac{1}{\sqrt{14}} \cdot\begin{pmatrix}1&3\\ 0&1\end{pmatrix}\otimes\begin{pmatrix}\cos \theta& -\sin \theta \\ \sin \theta & \cos \theta\end{pmatrix},\qquad A_2:=A_1^T,\]
where $\theta/\pi$ is irrational. Then $\A$ is strongly irreducible, $\dimaff \A=2$, and for every orthogonal projection $Q \in M_2(\R)$ we have $\dim_{\mathsf{aff}}^{I\otimes Q} \A \leq \log(625/4) / \log(125/2)\simeq 1.12966577\ldots$
\end{proposition}
\begin{figure}
    \centering
    \subfloat[Orthogonal projection of the attractor onto the first two co-ordinates.]{{\includegraphics[width=5.7cm]{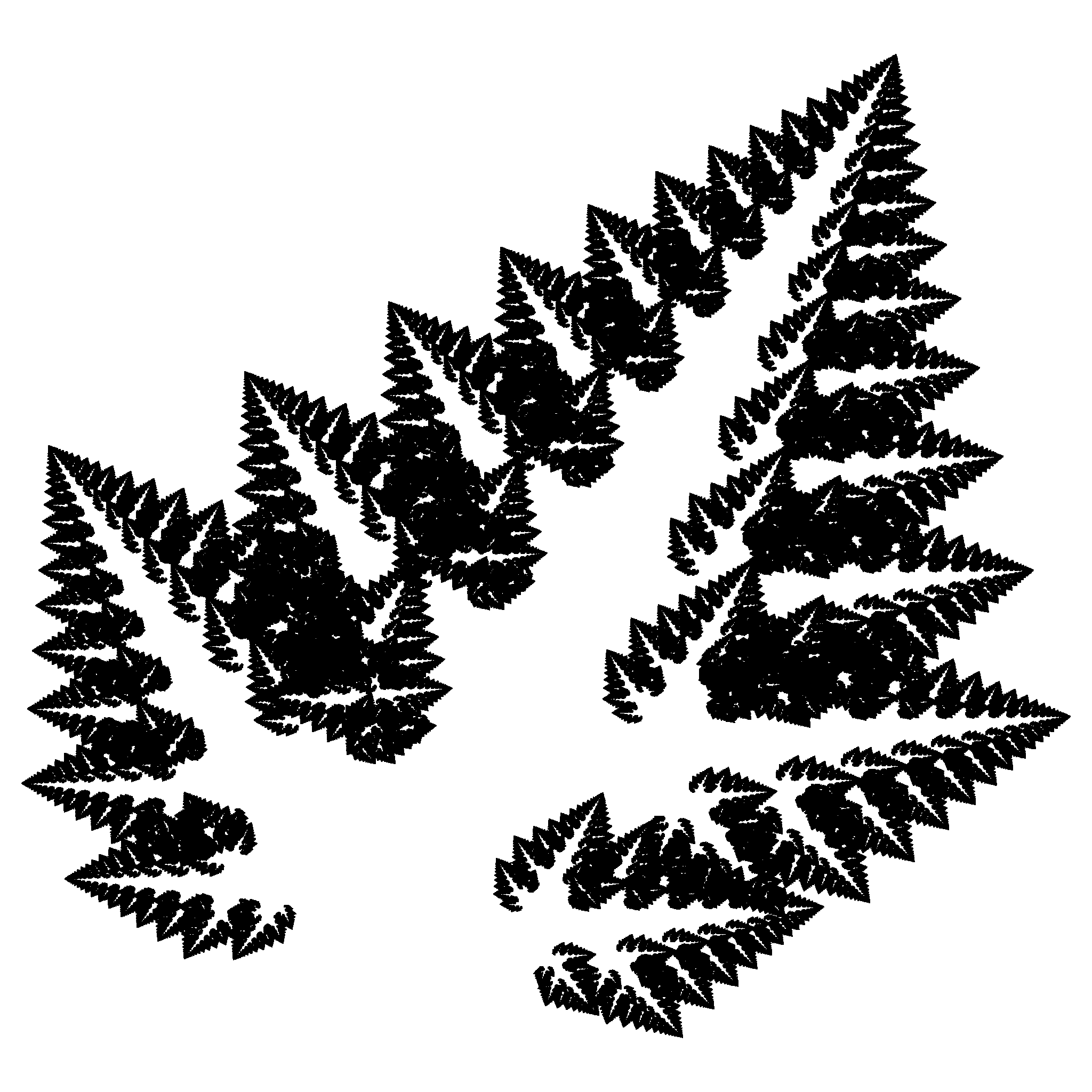}}}
    \quad
    \subfloat[Orthogonal projection of the attractor onto the first and third co-ordinates.]{{\includegraphics[width=5.7cm]{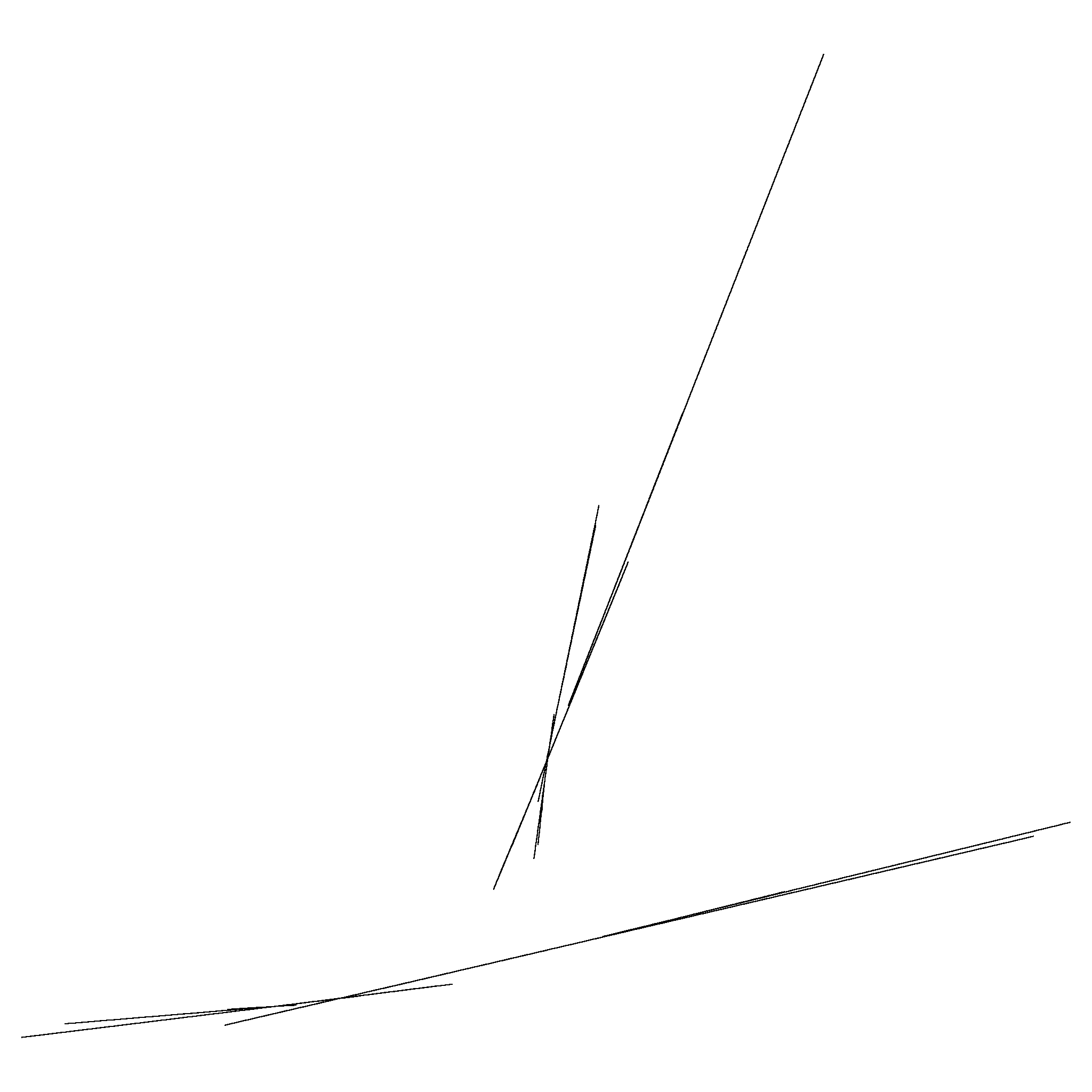}}}
    \caption{Two projections of the attractor of a system $(T_1, T_2)$ whose linearisation is the example $\A$ defined in Proposition \ref{pr:old} and whose additive parts are given by $v_1:=0$ and $v_2:=(1,0,1,0)^T$. The first projection corresponds to a linear map of the form $Q \otimes I$ and the second to a map of the form $I\otimes Q$; the latter results in a substantial apparent dimension drop.}
    \label{fi:gure}%
\end{figure}
The fact that this example has affinity dimension precisely $2$ makes it an interesting test case for the dimension problem for self-affine sets: by computing and examining projections of attractors of affine iterated function systems with this linearisation, one might hope to make a reasonable guess as to whether or not the attractor has dimension $2$. It was precisely this question which led, ultimately, to the results described in this article: in mid-2023 the author computed the projections displayed in Figure \ref{fi:gure}, observed the dimension drop phenomenon, and began the programme of research which led to the article \cite{MoSe25}. 
The existence of self-affine sets with tensor product structure and which admit exceptional projections as above was noted as a remark in \cite{MoSe25}, but a detailed construction and analysis was not published prior to the present work.

Before proceeding to the proof of Proposition \ref{pr:old} we establish a lemma which will support the precise calculation of the pressure. Variations of the following result have been applied across diverse works such as \cite{BlNe05,CrPaVu93,FeLoSh20,JuPr11,Mo19}. Here $\rho(B)$ denotes the spectral radius of the matrix $B$.
\begin{lemma}\label{le:plus}
Let $\B=(B_1,\ldots,B_N) \in M_d(\R)^N$ consist of non-negative matrices, and let $k \geq 1$. Then
\[\lim_{n \to \infty}\frac{1}{n}\log \sum_{|\iii|=n} \|B_\iii\|^k = \log\rho\left(\sum_{i=1}^N B_i^{\otimes k}\right).\]
\end{lemma}
\begin{proof}
Let $\threebar{\cdot}_1$ denote the entry-sum norm on $M_{d^k}(\R)$, i.e. define $\threebar{A}_1$ to be the sum of the absolute values of the matrix entries of $A$. Clearly $\threebar{A_1+A_2}_1= \threebar{A_1}_1+\threebar{A_2}_1$ when the matrices  $A_1,A_2 \in M_{d^k}(\R)$ are non-negative. Since every $B_i^{\otimes k}$ is non-negative,
\begin{align*}\lim_{n \to \infty}\frac{1}{n}\log \sum_{|\iii|=n} \|B_\iii\|^k &= \lim_{n \to \infty}\frac{1}{n}\log \sum_{|\iii|=n} \left\|B_\iii^{\otimes k}\right\|\\
&=\lim_{n \to \infty}\frac{1}{n}\log \sum_{|\iii|=n} \threebar{B_\iii^{\otimes k}}_1\\
& =\lim_{n \to \infty}\frac{1}{n}\log \threebar{\sum_{|\iii|=n} B_\iii^{\otimes k}}_1\\
& =\lim_{n \to \infty}\frac{1}{n}\log \threebar{\left(\sum_{i=1}^N B_i^{\otimes k}\right)^n}_1 = \log\rho\left(\sum_{i=1}^N B_i^{\otimes k}\right)\end{align*}
using Gelfand's formula.
\end{proof}

\begin{proof}[Proof of Proposition \ref{pr:old}]
Define
\[B_1:=\frac{1}{\sqrt{14}} \cdot\begin{pmatrix}1&3\\ 0&1\end{pmatrix}, \qquad B_1':=\begin{pmatrix}\cos \theta& -\sin \theta \\ \sin \theta & \cos \theta\end{pmatrix}, \quad B_2:=B_1^T,\quad B_2':=B_2^T \]
so that $A_i=B_i \otimes B_i'$ for $i=1,2$. Let us show that $\A$ is strongly irreducible. 
Suppose for a contradiction that $V_1,\ldots,V_k$ are nontrivial linear subspaces of $\R^4$ which are permuted by $\A$. Since $A_1$ and $A_2$ in particular induce some permutation of the set $\{V_1,\ldots,V_k\}$, by elementary group theory their powers $A_1^{k!}$ and $A_2^{k!}$ must induce the identity permutation of $\{V_1,\ldots,V_k\}$. Thus $A_1^{k!}$ and $A_2^{k!}$ share a nontrivial invariant subspace $V\subset \R^4$, say. Since $\theta/\pi $ is irrational both matrices have only non-real eigenvalues, so $\dim V$ cannot be odd and so by elimination must equal $2$. Such a subspace must be spanned by the real and imaginary parts of some complex eigenvector of $A_1^{k!}$, and must also be spanned by the real and imaginary parts of some complex eigenvector of $A_2^{k!}$. The space $V$ therefore has the form $\{u \otimes v \colon v\in \R^2\}$ for some common eigenvector $u \in \R^2$ of the matrices $B_1^{k!}$ and $B_2^{k!}$. The eigenvector $u$ is obviously also a common eigenvector of $B_1$ and $B_2$: but  no such common eigenvector exists. This  contradiction demonstrates that $\A$ is strongly irreducible as claimed. 

We now establish the various dimension estimates. For every word $\iii$ the two singular values of $B_\iii'$ are both $1$, so the four singular values $\sigma_i(B_\iii)\sigma_j(B_\iii')$ of $A_\iii$ are $\sigma_1(B_\iii)$, $\sigma_1(B_\iii)$, $\sigma_2(B_\iii)$ and $\sigma_2(B_\iii)$. Thus $\varphi^2(A_\iii)\equiv \|B_\iii\|^2$, so by Lemma \ref{le:plus},
\[P(\A,s)=\lim_{n \to \infty}\frac{1}{n}\log \sum_{|\iii|=n}\varphi^2(A_\iii) = \lim_{n \to \infty}\frac{1}{n}\log \sum_{|\iii|=n}\|B_\iii\|^2=\log\rho\left(B_1^{\otimes 2}+B_2^{\otimes 2}\right)\]
and an elementary calculation shows that the latter is zero. Thus $\dimaff \A=2$. Now let $Q \in M_d(\R)$ be an orthogonal projection. For every $\iii$ the singular values of $(I\otimes Q)(B_\iii \otimes B_\iii')$ are $\sigma_1(B_\iii)$, $\sigma_2(B_\iii)$, $0$ and $0$, so $\varphi^2((I\otimes Q)A_\iii)\equiv |\det B_\iii|$ and therefore 
\[P_{I\otimes Q}(\A,2)=\lim_{n \to \infty}\frac{1}{n}\log\sum_{|\iii|=n} \varphi^2((I\otimes Q)A_\iii) = \log\sum_{i=1}^2 |\det B_i| = \log \frac{1}{7}<0.\]
By similar calculations $\varphi^1((I\otimes Q)A_\iii) \equiv \|B_\iii\|$, so by Lemma \ref{le:plus}
\[P_{I\otimes Q}(\A,1)= \lim_{n \to \infty}\frac{1}{n}\log \sum_{|\iii|=n} \|B_\iii\| = \log \rho(B_1+B_2)=\log \frac{5}{\sqrt{14}}>0.\]
Since $s \mapsto P_{I\otimes Q}(\A,s)$ is convex on unit-length intervals with integer endpoints, interpolation between these two values demonstrates that
\[P_{I\otimes Q}\left(\A,\frac{\log 625-\log 4}{\log 175-\log 2}\right) \leq 0\]
and the claimed bound on $\dim_{\mathsf{aff}}^{I\otimes Q} \A$ follows.
\end{proof}

\subsection{A general example constructed via Kronecker products}
We complement the preceding concrete example with a more general construction as follows. 
\begin{theorem}\label{th:kronecker}
Let $\B:=(B_1,\ldots,B_N) \in \GL_{d_1}(\R)^N$ and $\B':=(B_1',\ldots,B_N')\in \GL_{d_2}(\R)^N$, and define $\A \in \GL_{d_1d_2}(\R)^N$ by $A_i:=B_i\otimes B_i'$ for every $i=1,\ldots,N$. Let $1 \leq k <d_1$ and $1 \leq \ell<d_2$, and suppose that $\A$ is contracting, that $\dimaff \A \in (k\ell, \min\{kd_2,d_1\ell\}]$, and that for some real number $\kappa>0$
\begin{equation}\label{eq:key} \frac{\sigma_1(B_\iii)}{\sigma_d(B_\iii)} <  e^{-\kappa |\iii|} \frac{\sigma_\ell(B_\iii')}{\sigma_{\ell+1}(B_\iii')}\end{equation}
for all long enough words $\iii$.  If $X$ is the attractor of an iterated function system with  linearisation $\A$, then for every orthogonal projection $Q \in M_{d_1}(\R)$ of rank $k$
\[\dimh (Q\otimes I)X < \dimaff \A.\]
\end{theorem}
If $B_i \in \O(d_1)$ for every $i$, and if  $\B'$ is $\ell$-dominated, then the hypothesis concerning \eqref{eq:key} is automatically satisfied.  Clearly, Theorem \ref{th:kronecker} implies the existence of self-affine sets for which every orthogonal projection of the form $Q\otimes I$ ---  where $\rank Q=k$, say --- is an exceptional projection.
\begin{proof}
Let $Q \in M_{d_1}(\R)$ have rank $k$, and define $s:=\dimaff \A \in (k\ell, \min\{kd_2,d_1\ell\}]$. Suppose that $\iii$ is long enough that \eqref{eq:key} holds. If $1 \leq i \leq d_1$ and $1 \leq j \leq \ell$, then using \eqref{eq:key}
 \[\sigma_i(B_\iii)\sigma_j(B_\iii')\geq \sigma_{d_1}(B_\iii)\sigma_\ell(B_\iii') > \sigma_{1}(B_\iii)\sigma_{\ell+1}(B_\iii')\geq \sigma_{i'}(B_\iii)\sigma_{j'}(B_\iii')\]
 for every $i',j'$ such that $1 \leq i' \leq d_1$ and $\ell+1 \leq j' \leq d_2$. This calculation demonstrates that the largest $d_1\ell$ singular values of $A_\iii$ must be the values $\sigma_i(B_\iii)\sigma_j(B_\iii')$ such that $1 \leq i \leq d_1$ and $1 \leq j \leq \ell$, in some unknown order. In particular it is definitely the case that
 \[\varphi^{k\ell}(A_\iii) \geq \prod_{\substack{1 \leq i \leq k \\ 1 \leq j \leq \ell}} \sigma_i(B_\iii)\sigma_j(B_\iii')\]
 and 
 \[\min_{k\ell+1 \leq i \leq d_1\ell} \sigma_i(A_\iii) \geq \sigma_{d_1}(B_\iii)\sigma_\ell(B_\iii'),\]
and consequently
\begin{equation}\label{eq:component-2} \left(\prod_{\substack{1 \leq i \leq k \\ 1 \leq j \leq \ell}} \sigma_i(B_\iii)\sigma_j(B_\iii')\right)\left( \sigma_d(B_\iii) \sigma_{\ell}(B_\iii')\right)^{s-k\ell}\leq \varphi^s(A_\iii)\end{equation}
 for all long enough words $\iii$. On the other hand, the singular values of $(Q\otimes I)A_\iii=(QB_\iii)\otimes B_\iii'$ are precisely the values $\sigma_i(QB_\iii)\sigma_j(B_\iii')$ for $i=1,\ldots,d_1$ and $j=1,\ldots,d_2$. Since  $\sigma_i(QB_\iii)=0$ when $k +1 \leq i \leq d_1$ there are exactly $kd_2$ nonzero singular values for $(QB_\iii)\otimes B_\iii'$. These $kd_2$ numbers, listed in decreasing order, are necessarily bounded above by the $kd_2$ values $\sigma_i(B_\iii)\sigma_j(B_\iii')$ such that $1 \leq i \leq k$ and $1 \leq j \leq d_2$, also listed in decreasing order. Now, if $1 \leq i \leq k$ and $1 \leq j \leq \ell$, then similarly to before
 \[\sigma_i(B_\iii)\sigma_j(B_\iii')\geq \sigma_{d_1}(B_\iii)\sigma_\ell(B_\iii') > \sigma_{1}(B_\iii)\sigma_{\ell+1}(B_\iii')\geq \sigma_{i'}(B_\iii)\sigma_{j'}(B_\iii')\]
 for any $i',j'$ such that $1 \leq i' \leq k$ and $\ell+1 \leq j' \leq d_2$. This demonstrates that the $k\ell$ largest expressions of the form $\sigma_i(B_\iii)\sigma_j(B_\iii')$ with  $1 \leq i \leq k$ and $1 \leq j \leq d_2$ must be those expressions such that $1 \leq i \leq k$ and $1 \leq j \leq \ell$, with all other values being bounded above by $\sigma_{1}(B_\iii)\sigma_{\ell+1}(B_\iii')$. Thus
 \[\varphi^{k\ell}((Q\otimes I)A_\iii) \leq \prod_{\substack{1 \leq i \leq k \\ 1 \leq j \leq \ell}} \sigma_i(B_\iii)\sigma_j(B_\iii')\]
 and
 \[\max_{k\ell+1 \leq i \leq kd_2} \sigma_i((Q\otimes I)A_\iii) \leq \sigma_1(B_\iii)\sigma_{\ell+1}(B_\iii')\]
 from which it follows that
 \begin{equation}\label{eq:component-1}\varphi^s((Q\otimes I)A_\iii) \leq \left(\prod_{\substack{1 \leq i \leq k \\ 1 \leq j \leq \ell}} \sigma_i(B_\iii)\sigma_j(B_\iii')\right)\left( \sigma_1(B_\iii) \sigma_{\ell+1}(B_\iii')\right)^{s-k\ell}\end{equation}
 for all long enough words $\iii$. Combining \eqref{eq:component-2} and \eqref{eq:component-1} with \eqref{eq:key} we find that for all large enough $n$,
\[\sum_{|\iii|=n} \varphi^s((Q\otimes I)A_\iii) \leq e^{-\kappa(s-k\ell)n} \sum_{|\iii|=n} \varphi^s(A_\iii)\]
so that
\[P_{Q\otimes I}(\A,s)\leq P(\A,s) - \kappa(s-k\ell) <P(\A,s)=0\]
 and applying Corollary \ref{co:cond} completes the proof of the theorem.
\end{proof}

\subsection{Examples arising from the split orthogonal group}
Our final class of examples is defined in terms of the \emph{split orthogonal group}.  For each $k \geq 1$, the split orthogonal group $\O(k,k)$ is defined to be the subgroup of $ \GL_{2k}(\R)$ consisting of all linear maps which preserve the bilinear form
\[f(u,v):=\sum_{i=1}^k u_iv_i - \sum_{i=k+1}^{2k} u_i v_i.\]
The subgroup $\SO(k,k)<\O(k,k)$ is defined to be the set of all elements of $\O(k,k)$ with unit determinant. An alternative view of this definition is as follows: if $J \in \GL_{2k}(\R)$ is defined to be the diagonal matrix whose first $k$ diagonal entries are $1$ and whose last $k$ diagonal entries are $-1$, then $B \in \O(k,k)$ if and only if $\langle Bu,JBv\rangle=\langle u,Jv\rangle$ for every $u,v \in \R^{2k}$, if and only if $B^TJB=J$. It is immediate from the latter characterisation that every $B \in \O(k,k)$ is conjugate via an isometry to its inverse transpose, hence satisfies $\sigma_i(B)=1/\sigma_{2k+1-i}(B)$ for every $i=1,\ldots,2k$ and has determinant $\pm1$. When treating the split orthogonal group it is often natural to view elements of $\R^{2k}$ as pairs $(u,v) \in \R^k \times \R^k$ and we will frequently adopt this viewpoint without comment.  

The group $\O(k,k)$ may be explicitly described in terms of its Cartan decomposition as follows:
\begin{proposition}[Cartan decomposition of the split orthogonal group]\label{pr:cartan}
Let $k \geq 1$, and define two subgroups of $\GL_{2k}(\R)$ by
\[\mathbf{A}=\left\{\begin{pmatrix}\cosh \theta_1 & \cdots & 0 & \sinh \theta_1 &\cdots &0\\ 
\vdots&\ddots &\vdots &\vdots& \ddots &\vdots\\
0& \cdots &\cosh \theta_k &0 &\cdots &\sinh \theta_k\\ \sinh \theta_1 & \cdots & 0 & \cosh \theta_1 &\cdots &0\\ 
\vdots&\ddots &\vdots &\vdots& \ddots &\vdots\\
0& \cdots &\sinh \theta_k &0 &\cdots &\cosh \theta_k
\end{pmatrix} \colon \theta_1,\ldots,\theta_k \in \R\right\}\]
and
\[\mathbf{K}:=\left\{\begin{pmatrix} U&0 \\ 0&V \end{pmatrix} \colon U,V \in \O(k)\right\}.\]
Then $\O(k,k)=\mathbf{KAK}$. 
\end{proposition}

Proposition \ref{pr:cartan} is a special case of the general Cartan decomposition for Lie groups and is most usually proved using Lie-theoretic machinery. For the convenience of a broad mathematical audience we present a quick proof via the properties of the standard singular value decomposition.
\begin{proof}
By direct calculations $\mathbf{K}$ and $\mathbf{A}$ are subgroups of $\O(k,k)$, so $\mathbf{KAK}\subseteq \O(k,k)$. Let us demonstrate the reverse inclusion. Given $B \in \O(k,k)$, write $B=\left(\begin{smallmatrix} A_{11} &A_{12}\\ A_{21}&A_{22}\end{smallmatrix}\right)$ where each $A_{ij}$ is a $k \times k$ block matrix. For $i=1,2$ let $C_i$ denote the diagonal matrix whose entries are the singular values of $A_{ii}$ in decreasing order. Appealing to the existence of singular value decompositions for $A_{11}, A_{22} \in M_k(\R)$ we may write
\[B=\begin{pmatrix} U_1 & 0 \\ 0&U_2 \end{pmatrix} \begin{pmatrix} C_1 &X_2 \\ X_1 & C_2\end{pmatrix} \begin{pmatrix} V_1 & 0 \\ 0&V_2 \end{pmatrix}\]
with $U_1, U_2, V_1, V_2 \in \O(k)$ and $X_1,X_2 \in M_k(\R)$.  Denote the middle one of the above three block matrices by $A \in \O(k,k)$. The equation $A^TJA=J$ implies the equations
\[C_1^2-X_1^TX_1=C_2^2-X_2^TX_2=I,\qquad C_1X_2=X_1^TC_2.\]
The first identity implies that the singular values of $C_1$ and of $C_2$ are all at least $1$, so write these as $\cosh \theta_{i,1} \geq \cdots \geq \cosh \theta_{i,k}$ where each $\theta_{i,j}$ is non-negative. It follows directly that $X_i^TX_i=C_i^2-I=S_i^2$ for $i=1,2$, where each $S_i$ is the diagonal matrix with entries $\sinh \theta_{i,1}\geq \cdots\geq \sinh \theta_{i,k} \geq 0$. By \cite[Theorem 7.3.11]{HoJo13} this implies that $X_i=O_iS_i$ for $i=1,2$ where each $O_i$ is orthogonal. The equation $C_1X_2=X_1^TC_2$ now implies $O_2S_2C_2^{-1}=C_1^{-1}S_1O_1^T$ which by equating singular values yields $\tanh \theta_{1,j}=\tanh \theta_{2,j}$ for every $j=1,\ldots,k$. Since $\tanh$ is injective we conclude that $\theta_{1,j}\equiv \theta_{2,j}$ and therefore $C_1=C_2=C$ and $S_1=S_2=S$, say.  Thus
\begin{equation}\label{eq:O}B=\begin{pmatrix} U_1 & 0 \\ 0&U_2 \end{pmatrix} \begin{pmatrix} C &O_1S \\ O_2S & C\end{pmatrix} \begin{pmatrix} V_1 & 0 \\ 0&V_2 \end{pmatrix}\end{equation}
where $O_2SC^{-1}O_1=SC^{-1}$. By \cite[Theorem 2.6.5]{HoJo13} the latter is only possible if each $O_i$ preserves every eigenspace of $SC^{-1}$, which implies preserving every eigenspace of $S$ and of $C$ since all three matrices have the same eigenspaces. Without loss of generality we may modify $O_1$ and $O_2$ so as to coincide with the identity on $\ker S$ while remaining unchanged on the other eigenspaces; clearly this does not affect the validity of \eqref{eq:O} nor the orthogonality of $O_1$ and $O_2$. Each $O_i$ preserves every eigenspace of $S$ and $C$ and consequently commutes with $C$ and $S$. From $O_2SC^{-1}O_1=SC^{-1}$ we deduce $O_2=O_1^{-1}$ and hence
$$B=\begin{pmatrix} U_1 & 0 \\ 0&U_2 \end{pmatrix} \begin{pmatrix} I & 0 \\ 0&O_1^{-1} \end{pmatrix} \begin{pmatrix} C &S \\ S & C\end{pmatrix}\begin{pmatrix} I& 0 \\ 0&O_1 \end{pmatrix} \begin{pmatrix} V_1 & 0 \\ 0&V_2 \end{pmatrix}\in \mathbf{KAK}$$ proving the proposition.\end{proof}

We will call a subspace $W\leq \R^{2k}$ \emph{isotropic} if $f(u,v)=0$ for every $u,v \in W$. In this article we will be exclusively concerned with isotropic subspaces of dimension $k$, and we will see shortly that these are precisely the subspaces of the form $\{(u,Ou) \in \R^{2k} \colon u \in \R^k\}$ where $O \in \O(k)$. We will call such a subspace \emph{positively oriented} if $\det O=1$, and \emph{negatively oriented} otherwise. We let $\R^* \O(k,k)$ (respectively $\R^*\SO(k,k)$) denote respectively the group of nonzero scalar multiples of elements of $\O(k,k)$ (respectively $\SO(k,k)$).
The main result of this section is the following:
\begin{theorem}\label{th:split-example}
Let $\A=(A_1,\ldots,A_N) \in (\R^* \SO(k,k))^N$ be contracting and $k$-dominated, with $\dimaff \A \in (k-1,k]$. Then either $\dimh QX<\dimaff \A$
 for every orthogonal projection $Q \in M_{2k}(\R)$ whose image is a positively oriented isotropic space of rank $k$, 
or the same inequality holds for every projection whose image is a \emph{negatively oriented} isotropic space of rank $k$.
\end{theorem}
 \begin{figure}
    \centering
    \subfloat{{\includegraphics[width=5.7cm]{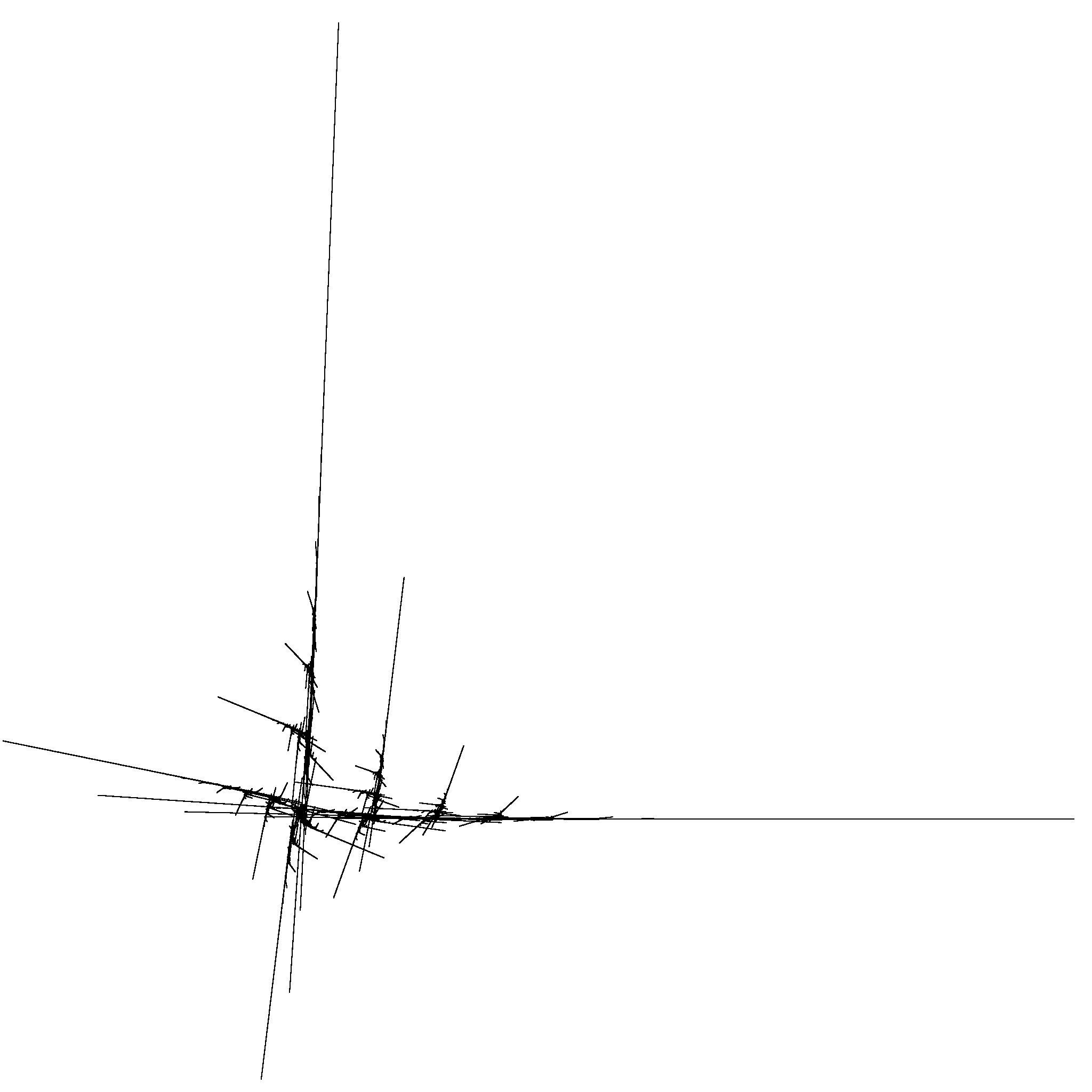}}}
\quad
    \subfloat{{\includegraphics[width=5.7cm]{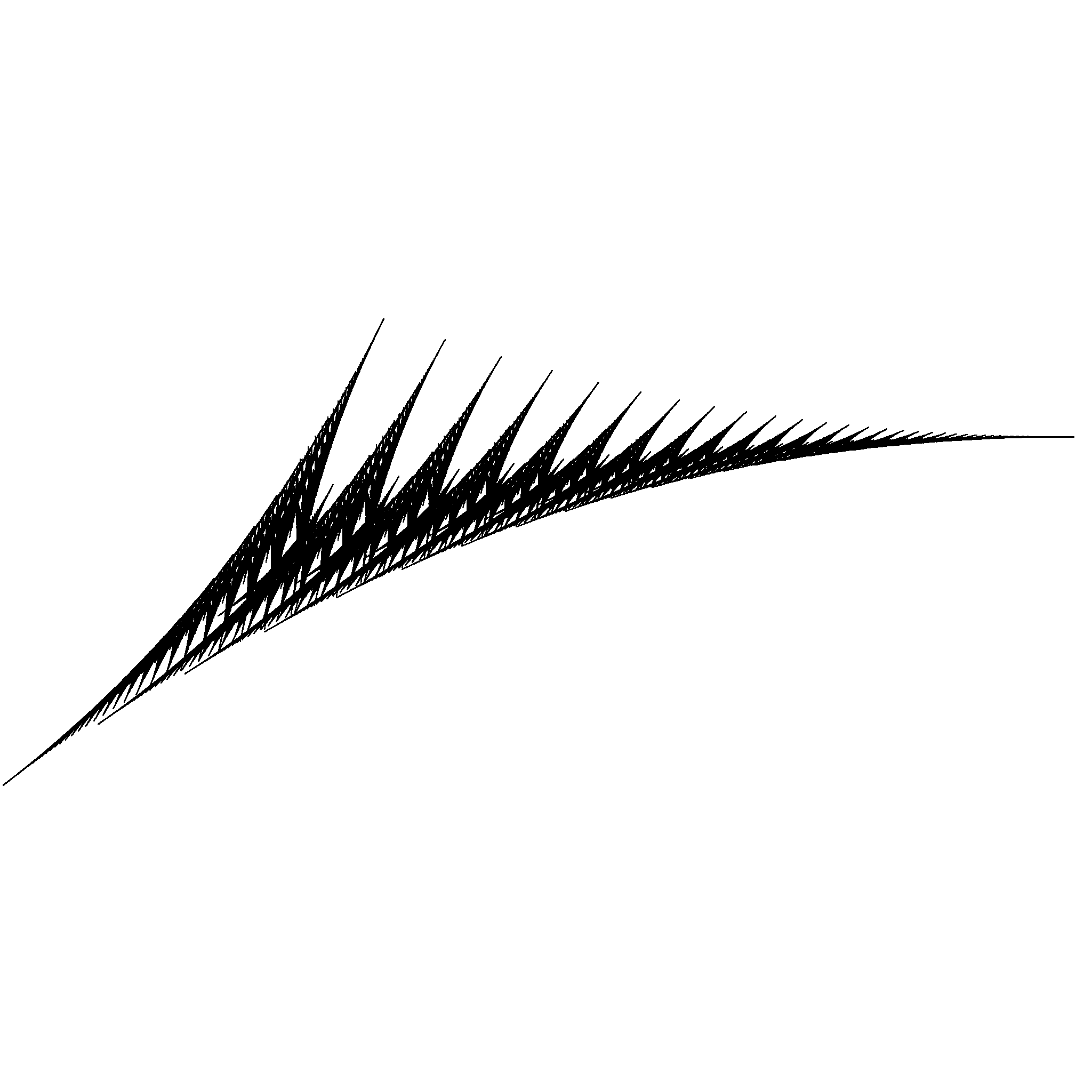}}}
    \caption{Two projections of the attractor of the system $(T_1, T_2)$ defined in \eqref{eq:so22-example} onto isotropic planes of differing orientations.}
    \label{fi:gure-two}%
\end{figure}
It follows that if $\dimh X=\dimaff \A$ then the set of exceptional projections of $X$ either includes all orthogonal projections whose kernel has the form $\{(u,Ou) \colon u \in \R^k\}$
where $\det O=1$, or includes all orthogonal projections whose kernel has the same form with $\det O=-1$. If $\A$ is strongly contracting then the existence of attractors satisfying $\dimh X=\dimaff \A$ is guaranteed by Falconer's theorem. As was the case with Theorem \ref{th:kronecker}, it is not difficult to construct cases of Theorem \ref{th:split-example} in which $\A$ acts strongly irreducibly on $\R^{2k}$.

Similarly to the tensor product constructions of the previous two subsections, the existence of exceptional projections of self-affine sets with linearisations in $\R^*\SO(k,k)$ was stated as a remark in \cite{MoSe25}, but no detailed construction has been given prior to the present work.

Figure \ref{fi:gure-two} presents two projections of the attractor of an iterated function system $(T_1,T_2)$ given by 
\begin{align}\label{eq:so22-example}T_1\begin{pmatrix}w\\x\\y\\z\end{pmatrix}&:=\frac{1}{2}\begin{pmatrix}\cosh \frac{3}{5}&0&\sinh\frac{3}{5}&0\\
0&\cosh\frac{1}{2} &0&\sinh\frac{1}{2}\\
\sinh\frac{3}{5} &0&\cosh\frac{3}{5}&0\\
0&\sinh\frac{1}{2} &0&\cosh\frac{1}{2}
\end{pmatrix}\begin{pmatrix}w\\x\\y\\z\end{pmatrix}+\begin{pmatrix}1\\0\\1\\0\end{pmatrix}\\\nonumber
T_2\begin{pmatrix}w\\x\\y\\z\end{pmatrix}&:=\frac{2}{3}\begin{pmatrix}\cos \sqrt{3}&-\sin \sqrt{3}&0&0\\
\sin \sqrt{3}&\cos\sqrt{3}&0&0\\
0&0&\cos 5&-\sin 5\\
0&0&\sin 5&\cos5
\end{pmatrix}\begin{pmatrix}w\\x\\y\\z\end{pmatrix}+\begin{pmatrix}0\\-1\\0\\-1\end{pmatrix}.
\end{align}
Although this example is not $2$-dominated, a dimension drop for certain projections is clearly visible empirically.

The underlying principle of Theorem \ref{th:split-example} is fundamentally the same as that of the previous examples: certain spaces spanned by singular vectors of products $A_\iii$ will be shown to necessarily intersect the kernels of certain projections $Q$. In this example, however, these intersections rely on these spaces being $k$-dimensional isotropic spaces with particular orientations. We therefore begin the process of proof by  demonstrating some general features of the geometry of isotropic subspaces and of their behaviour under the action of the groups $\O(k,k)$ and $\SO(k,k)$. We first verify the description of these subspaces stated previously:
\begin{lemma}\label{le:isotropic-description}
Let $k \geq 1$ and let $W \leq \R^{2k}$ have dimension $k$. Then $W$ is isotropic if and only if there exists $O \in \O(k)$ such that $W=\{(u,Ou) \colon u \in \R^k\}$.
\end{lemma}
\begin{proof}
Suppose $W <\R^{2k}$ is an isotropic subspace of dimension $k$, and define $P \colon \R^{2k} \to \R^{2k}$ by orthogonal projection onto the first $k$ co-ordinates, as $(u,v) \mapsto (u,0)$. Clearly $\|Pw\|^2 - \|(I-P)w\|^2=f(w,w)=0$ for every $w \in W$, and since also $\|Pw\|^2+\|(I-P)w\|^2=\|w\|^2$ it follows that $P$ defines an injection $W\to PW$, hence a bijection. By dimension considerations $PW$ must be precisely equal to $\{(u,0) \colon u \in \R^k\}$. The inverse transformation $PW \to W$ is linear and so takes each $(u,0)$ to some $(u,Lu)$ for some consistent $L \in M_k(\R)$, whence necessarily $W=\{(u,Lu) \colon u \in \R^k\}$. Since $f((u,Lu),(v,Lv))=\langle u,v\rangle-\langle Lu,Lv\rangle=0$ for all $u,v \in \R^k$ the linear map $L$ is an isometry, and this shows that every $k$-dimensional isotropic subspace has the required form. It is conversely obvious that every subspace of the form stated is isotropic.\end{proof}
The intersection properties of isotropic subspaces will form a pivotal step in the proof of Theorem \ref{th:split-example}:
\begin{lemma}\label{le:isotropic-intersection}
Let $k \geq 1$ and let $W_1, W_2<\R^{2k}$ be isotropic subspaces of dimension $k$ having the form $W_i=\{(u,O_iu) \colon u \in \R^k\}$ for each $i=1,2$, where $O_i \in \O(k)$. Then $\dim (W_1 \cap W_2)$ is odd if and only if $(\det O_1)(\det O_2) = (-1)^{k+1}$.
\end{lemma}
\begin{proof}
We have $(u,v) \in W_1 \cap W_2$ if and only if $O_1u=O_2u$ if and only if $u$ is fixed by $O_1^{-1}O_2$, so $\dim (W_1 \cap W_2)$ is precisely the multiplicity of $1$ as an eigenvalue of $O_1^{-1}O_2$. Since  $O_1^{-1}O_2$ is an orthogonal matrix, its eigenvalues are $1$ and $-1$ with some unknown multiplicities, together with a (possibly empty) set of conjugate pairs of complex units. If $\det(O_1^{-1}O_2)=1$, therefore,  the multiplicity of $-1$ must be even: in this case $\dim (W_1 \cap W_2)$ is odd if and only if $k$ is also odd. Similarly, if $\det(O_1^{-1}O_2)=-1$ then the multiplicity of $-1$ is odd, in which case $\dim (W_1 \cap W_2)$ is odd if and only if $k$ is even.
\end{proof}
\begin{lemma}\label{le:isotropic-switching}
If $W <\R^{2k}$ is an isotropic subspace of dimension $k$, and $B \in \O(k,k)$, then $BW$ is isotropic, and its orientation matches that of $W$ if and only if $\det B=1$.
\end{lemma}
\begin{proof}
That the subspace $BW$ is isotropic follows directly from the definitions. Using Lemma \ref{le:isotropic-description}, write $W=\{(u,Ou) \colon u \in \R^k\}$ where $O \in \O(k)$.  To prove the remaining clause, we note that by Proposition \ref{pr:cartan} it suffices for us to check that every $B \in \mathbf{A}$ preserves the orientation of $W$, and that $B=\left(\begin{smallmatrix}U&0\\ 0&V\end{smallmatrix}\right) \in \mathbf{K}$ preserves the orientation of $W$ if and only if $(\det U)(\det V)=1$. In the latter case we note that $BW=\{(u,VOU^{-1}u) \colon u \in \R^k\}$, and clearly $\det (VO U^{-1})=(\det U)(\det V)(\det O)$ which makes the conclusion obvious. In the former case, using the definition of $\mathbf{A}$ we may write $B=\left(\begin{smallmatrix}C&S\\S&C\end{smallmatrix}\right)$ where $C$ and $S$ are commuting symmetric matrices satisfying $C^2-S^2=I$. We have $BW=\{(u,(CO+S)(C+SO)^{-1}u) \colon u \in \R^k\}$, where by Lemma \ref{le:isotropic-description} (or a laborious calculation using the properties of $C$ and $S$) the matrix $(CO+S)(C+SO)^{-1}$ is necessarily orthogonal. To see that $\det ((CO+S)(C+SO)^{-1})=\det O \in \{\pm 1\}$ we note that this result is obvious in the case where $B$ is the identity element of $\mathbf{A}$, and the general case follows since $\mathbf{A}$ is  connected.
\end{proof}
Our final lemma relates singular values of elements of $\R^*\O(k,k)$ to isotropic spaces:
\begin{lemma}\label{le:singular-spaces}
For every $B \in \R^*\O(k,k)$ such that $\sigma_{k+1}(B)<\sigma_k(B)$, let $\mathcal{V}(B)$ denote the span of the eigenvectors of $B^TB$ corresponding to the eigenvalues $\sigma_{k+1}(B)^2,\ldots,\sigma_{2k}(B)^2$. Then $\mathcal{V}(B)$ is isotropic. If $B_1 \in \R^*\SO(k,k)$,  $B_2 \in \R^*\O(k,k)$ satisfy $(\|B_1\|\cdot \|B_1^{-1}\|) \sigma_{k+1}(B_2)<\sigma_k(B_2)$, then the subspaces $\mathcal{V}(B_1B_2)$, $\mathcal{V}(B_2)$ and $\mathcal{V}(B_2B_1)$ are well-defined and have the same orientation.
\end{lemma}
\begin{proof}
Clearly $\mathcal{V}(B)=\mathcal{V}(\lambda B)$ for any $\lambda \in \R^*$ and any $B \in \O(k,k)$ such that $\sigma_{k+1}(B)<\sigma_k(B)$, so we assume without loss of generality that $B, B_1, B_2$ belong to $\O(k,k)$.  If $B \in \O(k,k)$ satisfies $\sigma_{k+1}(B)<\sigma_k(B)$, then since $B$ is conjugated by an isometry to its inverse transpose, $\sigma_{k+1}(B)<1<\sigma_k(B)$. It follows that in this case $\mathcal{V}(B)$ is spanned precisely by those eigenvectors of $B^TB$ corresponding to eigenvalues strictly less than $1$. If $v_1,v_2 \in \mathcal{V}(B)$ are eigenvectors with eigenvalues $\lambda_1, \lambda_2$, say, then $f(v_1,v_2)=f(B^TB v_1, B^TB v_2) = \lambda_1 \lambda_2 f(v_1,v_2)$ implies $f(v_1,v_2)=0$ since $\lambda_1\lambda_2 \neq 1$. It follows that $\mathcal{V}(B)$ is isotropic. It is an easy consequence of the identity $B^TB = B^{-1}(BB^T) B$ and the definition of $\mathcal{V}$ that $B\mathcal{V}(B^T)=\mathcal{V}(B)$. 
Now let $B_1 \in \SO(k,k)$ and $B_2 \in \O(k,k)$ satisfy $(\|B_1\|\cdot \|B_1^{-1}\|) \sigma_{k+1}(B_2)<\sigma_k(B_2)$. By elementary calculations $\sigma_{k+1}(B)<\sigma_k(B)$ where $B$ is any of $B_1 B_2$, $B_2$ or $B_2B_1$. We claim that $\mathcal{V}(B_1B_2)$ and $\mathcal{V}(B_2)$ have the same orientation. Write these spaces respectively as $\{(u_1,O_1 u_1) \colon u_1 \in \R^k\}$ and $\{(u_2, O_2u_2) \colon u_2 \in \R^k\}$, and let $v=(u_1, O_1u_1) \in \mathcal{V}(B_1B_2)$ be a unit vector, in which case $\|u_1\|=1/\sqrt{2}$. Write $v=w_1+w_2$ with $w_1 \in \mathcal{V}(B_2)$ and $w_2 \in \mathcal{V}(B_2)^\perp$. We compute that
\begin{align*}\sigma_k(B_2) \|w_2\|\leq \|B_2 w_2\|\leq \|B_2 v\| &\leq \|B_1^{-1}\| \|B_1 B_2 v\|\\
& \leq \|B_1^{-1}\| \sigma_{k+1}(B_1 B_2)\\
&\leq   (\|B_1^{-1}\| \cdot\|B_1\|) \sigma_{k+1}(B_2) < \sigma_k(B_2)\end{align*}
so that $\|w_2\|<1$. Writing $w_1=(u_2,O_2u_2) \in V(B_2)$, we deduce that
\[\|u_1-u_2\|^2 + \|O_1u_1 - O_2u_2\|^2=\|v-w_1\|^2=\|w_2\|^2<1\]
whence
\[\|u_1-O_2^{-1}O_1u_1\| \leq \|u_1-u_2\| + \|u_2-O_2^{-1}O_1u_1\|<\sqrt{2}= 2\|u_1\|.\]
Since $u_1 \in \R^k$ was an arbitrary vector of norm $1/\sqrt{2}$, it follows that $-1$ is not an eigenvalue of $O_2^{-1}O_1$ and therefore $\det O_2^{-1}O_1=1$. The claim follows. To see that $\mathcal{V}(B_2B_1)$ also has the same orientation as $\mathcal{V}(B_2)$ we note that this is equivalent to $B_1^TB_2^T\mathcal{V}(B_1^TB_2^T)$ having the same orientation as $B_2^T\mathcal{V}(B_2^T)$, which by Lemma \ref{le:isotropic-switching} is equivalent to  $\mathcal{V}(B_1^TB_2^T)$ having the same orientation as $\mathcal{V}(B_2^T)$; but this follows directly by the applying the preceding claim to $B_1^T$, $B_2^T$ in place of $B_1, B_2$.\end{proof}
\begin{proof}[Proof of Theorem \ref{th:split-example}]
Choose $n_0 \geq 1$ large enough that $(\|A_j^{-1}\| \cdot \|A_j\|) \sigma_{k+1}(A_\iii)<\sigma_k(A_\iii)$ for all words $\iii$ of length at least $n_0$ and all $j =1,\ldots,N$, and define $\mathcal{V}(A_\iii^T)$ as in Lemma \ref{le:singular-spaces} for every such word $\iii$. We claim that the orientation of the isotropic space $\mathcal{V}(A_\iii^T)$ is consistent across all words $\iii$ of length at least $n_0$. To see this, let $|\iii|, |\jjj| \geq n_0$ and observe that by repeated application of Lemma \ref{le:singular-spaces} the subspaces $\mathcal{V}(A_\iii^T)$, $\mathcal{V}(A_{j_1}^TA_\iii^T )$, $\mathcal{V}(A_{j_2}^TA_{j_1}^TA_\iii^T)$,\ldots, $\mathcal{V}(A_\jjj^T A_\iii^T)$ share a common orientation. Removing symbols from $\iii$ one at a time in a similar manner we see that $\mathcal{V}(A_\jjj^TA_\iii^T), \mathcal{V}(A_\jjj^T A_{i_{|\iii|}}^T\cdots A_{i_2}^T),\ldots,\mathcal{V}(A_\jjj^TA_{i_{|\iii|}}^T)$ all have the same orientation as $\mathcal{V}(A_\jjj^T)$. Thus $\mathcal{V}(A_\iii^T)$ and $\mathcal{V}( A_\jjj^T)$ have the same orientation, which proves the claim. Now let $Q$ be an arbitrary orthogonal projection of rank $k$ whose image is an isotropic space matching the shared orientation of every $\mathcal{V}(A_\iii^T)$ if $k$ is odd, or opposite to the shared orientation of every $\mathcal{V}(A_\iii^T)$ if $k$ is even. By Lemma \ref{le:isotropic-intersection}, the image of $Q$ nontrivially intersects every subspace of the form $\mathcal{V}(A_\iii^T)$ where $|\iii|\geq n_0$. 

By Corollary \ref{co:cond}, to prove the theorem it is enough to show that $P_Q(\A,s)<P(\A,s)$, where $s:=\dimaff \A \in (k-1,k]$.  Let $|\iii|\geq n_0$. The singular values of $QA_\iii$ are precisely the square roots of the singular values of $QA_\iii A_\iii^TQ$, of which precisely $k$ are zero. By consideration of the singular value decomposition of the latter, $\sigma_k(QA_\iii)^2$ is precisely the infimal value of $\|QA_\iii A_\iii^T v\|$ over unit vectors $v$ in the image of $Q$. Choosing any  unit vector $v\in \im Q \cap \mathcal{V}(A_\iii^T)$ we find that
\[\sigma_{k}(QA_\iii)^2 =\sigma_{k}\left(QA_\iii A_\iii^TQ\right) \leq \left\|A_\iii A_\iii^T v\right\| \leq \sigma_{k+1}(A_\iii)^2\]
by the definition of $\mathcal{V}(A_\iii^T)$, for every word $\iii$ of length at least $n_0$. Since $\sigma_\ell(QB) \leq \sigma_\ell(B)$ for all $B \in M_{2k}(\R)$ and $\ell \leq k$, we conclude that
\begin{align*}P_Q(\A,s)&=\lim_{n \to \infty}\frac{1}{n}\sum_{|\iii|=n} \varphi^s(QA_\iii)\\ &\leq \lim_{n \to \infty}\frac{1}{n}\sum_{|\iii|=n}\left(\frac{\sigma_{k+1}(A_\iii)}{\sigma_k(A_\iii)}\right)^{s-(k-1)}\varphi^s(A_\iii)<P(\A,s)=0\end{align*}
using $k$-domination.
\end{proof}
\section{Acknowledgements}
Theorems \ref{th:kronecker} and \ref{th:split-example} originate in an earlier manuscript of the article \cite{MoSe25} which was some 30 pages longer than the version submitted for publication and distributed on arXiv. Using the results of \cite{MoSe25} those two theorems may be obtained much more concisely -- though perhaps less transparently -- than by the method presented here. The author thanks \cagri Sert for agreeing to the repurposing of this material for use in the present work.
\bibliographystyle{acm}
\bibliography{badminton}

\begin{thebibliography}{10}

\bibitem{BaHoRa19}
Bal\'{a}zs B\'{a}r\'{a}ny, Michael Hochman, and Ariel Rapaport.
\newblock Hausdorff dimension of planar self-affine sets and measures.
\newblock {\em Invent. Math.}, 216(3):601--659, 2019.

\bibitem{BaKe24}
Bal\'azs B\'ar\'any and Antti K\"aenm\"aki.
\newblock Dimension of planar non-conformal attractors with triangular
  derivative matrices.
\newblock {\em Comm. Math. Phys.}, 405(11):Paper No. 256, 31, 2024.

\bibitem{BaSiSo23}
Bal\'{a}zs B\'{a}r\'{a}ny, K\'{a}roly Simon, and Boris Solomyak.
\newblock {\em Self-similar and self-affine sets and measures}, volume 276 of
  {\em Mathematical Surveys and Monographs}.
\newblock American Mathematical Society, Providence, RI, 2023.

\bibitem{Be84}
Tim Bedford.
\newblock {\em Crinkly curves, {M}arkov partitions and dimension}.
\newblock 1984.
\newblock Thesis (Ph.D.)--The University of Warwick.

\bibitem{BlNe05}
Vincent~D. Blondel and Yurii Nesterov.
\newblock Computationally efficient approximations of the joint spectral
  radius.
\newblock {\em SIAM J. Matrix Anal. Appl.}, 27(1):256--272, 2005.

\bibitem{BoGo09}
Jairo Bochi and Nicolas Gourmelon.
\newblock Some characterizations of domination.
\newblock {\em Math. Z.}, 263(1):221--231, 2009.

\bibitem{CrPaVu93}
Andrea Crisanti, Giovanni Paladin, and Angelo Vulpiani.
\newblock {\em Products of random matrices in statistical physics}, volume 104
  of {\em Springer Series in Solid-State Sciences}.
\newblock Springer-Verlag, Berlin, 1993.

\bibitem{DaSi17}
Tushar Das and David Simmons.
\newblock The {H}ausdorff and dynamical dimensions of self-affine sponges: a
  dimension gap result.
\newblock {\em Invent. Math.}, 210(1):85--134, 2017.

\bibitem{Ed92}
Gerald~A. Edgar.
\newblock Fractal dimension of self-affine sets: some examples. \emph{Measure
  Theory (Oberwolfach, 1990)}.
\newblock {\em Rend. Circ. Mat. Palermo (2) Suppl.}, (28):341--358, 1992.

\bibitem{Fa88}
Kenneth~J. Falconer.
\newblock The {H}ausdorff dimension of self-affine fractals.
\newblock {\em Math. Proc. Cambridge Philos. Soc.}, 103(2):339--350, 1988.

\bibitem{FaOr14}
Katrin F\"{a}ssler and Tuomas Orponen.
\newblock On restricted families of projections in {$\mathbb R^3$}.
\newblock {\em Proc. Lond. Math. Soc. (3)}, 109(2):353--381, 2014.

\bibitem{Fe23}
De-Jun Feng.
\newblock Dimension of invariant measures for affine iterated function systems.
\newblock {\em Duke Math. J.}, 172(4):701--774, 2023.

\bibitem{FeLoSh20}
De-Jun Feng, Chiu-Hong Lo, and Shuang Shen.
\newblock Uniformity of {L}yapunov exponents for non-invertible matrices.
\newblock {\em Ergodic Theory Dynam. Systems}, 40(9):2399--2433, 2020.

\bibitem{FeXi25}
De-Jun Feng and Yu-Hao Xie.
\newblock Dimensions of orthogonal projections of typical self-affine sets and
  measures.
\newblock Preprint: arXiv:2502.04000, 2025.

\bibitem{Fr12}
Jonathan~M. Fraser.
\newblock On the packing dimension of box-like self-affine sets in the plane.
\newblock {\em Nonlinearity}, 25(7):2075--2092, 2012.

\bibitem{FrJu21}
Jonathan~M. Fraser and Natalia Jurga.
\newblock The box dimensions of exceptional self-affine sets in
  {$\mathbb{R}^3$}.
\newblock {\em Adv. Math.}, 385:Paper No. 107734, 32, 2021.

\bibitem{GaGuGuHaMaWa22}
Shengwen Gan, Shaoming Guo, Larry Guth, Terence L.~J. Harris, Dominique
  Maldague, and Hong Wang.
\newblock On restricted projections to planes in $\mathbb{R}^3$.
\newblock {\em Amer. J. Math.}
\newblock To appear. Preprint: arXiv:2207.13844.

\bibitem{Ho14}
Michael Hochman.
\newblock On self-similar sets with overlaps and inverse theorems for entropy.
\newblock {\em Ann. of Math. (2)}, 180(2):773--822, 2014.

\bibitem{Ho20}
Michael Hochman.
\newblock On self-similar sets with overlaps and inverse theorems for entropy
  in $\mathbb{R}^d$.
\newblock {\em Mem. Amer. Math. Soc.}, 265(1287), 2020.

\bibitem{HoJo94}
Roger~A. Horn and Charles~R. Johnson.
\newblock {\em Topics in matrix analysis}.
\newblock Cambridge University Press, Cambridge, 1994.
\newblock Corrected reprint of the 1991 original.

\bibitem{HoJo13}
Roger~A. Horn and Charles~R. Johnson.
\newblock {\em Matrix analysis}.
\newblock Cambridge University Press, Cambridge, second edition, 2013.

\bibitem{JuPr11}
Rapha\"el~M. Jungers and Vladimir~Yu. Protasov.
\newblock Fast methods for computing the {$p$}-radius of matrices.
\newblock {\em SIAM J. Sci. Comput.}, 33(3):1246--1266, 2011.

\bibitem{Ka04}
Antti K\"{a}enm\"{a}ki.
\newblock On natural invariant measures on generalised iterated function
  systems.
\newblock {\em Ann. Acad. Sci. Fenn. Math.}, 29(2):419--458, 2004.

\bibitem{KaOrVe25}
Antti K\"aenm\"aki, Tuomas Orponen, and Laura Venieri.
\newblock A {M}arstrand-type restricted projection theorem in $\mathbb{R}^3$.
\newblock {\em Amer. J. Math.}, 147(1):81--123, 2025.

\bibitem{Ma54}
John~M. Marstrand.
\newblock Some fundamental geometrical properties of plane sets of fractional
  dimensions.
\newblock {\em Proc. London Math. Soc. (3)}, 4:257--302, 1954.

\bibitem{Ma75}
Pertti Mattila.
\newblock Hausdorff dimension, orthogonal projections and intersections with
  planes.
\newblock {\em Ann. Acad. Sci. Fenn. Ser. A I Math.}, (no. 2,):227--244, 1975.

\bibitem{Mc84}
Curt McMullen.
\newblock The {H}ausdorff dimension of general {S}ierpi\'nski carpets.
\newblock {\em Nagoya Math. J.}, 96:1--9, 1984.

\bibitem{Mo19}
Ian~D. Morris.
\newblock An explicit formula for the pressure of box-like affine iterated
  function systems.
\newblock {\em J. Fractal Geom.}, 6(2):127--141, 2019.

\bibitem{Mo23}
Ian~D. Morris.
\newblock On affine iterated function systems which robustly admit an invariant
  affine subspace.
\newblock {\em Proc. Amer. Math. Soc.}, 151(1):101--112, 2023.

\bibitem{MoSe23}
Ian~D. Morris and Cagri Sert.
\newblock A variational principle relating self-affine measures to self-affine
  sets.
\newblock Preprint: arXiv:2303.03437.

\bibitem{MoSe25}
Ian~D. Morris and Cagri Sert.
\newblock Projections of self-affine fractals.
\newblock Preprint: arXiv:2502.04001, 2025.

\bibitem{PrYaZa22}
Malabika Pramanik, Tongou Yang, and Joshua Zahl.
\newblock A {F}urstenberg-type problem for circles, and a {K}aufman-type
  restricted projection theorem in $\mathbb{R}^3$.
\newblock {\em Amer. J. Math.}
\newblock To appear. Preprint: arXiv:2207.02259.

\bibitem{Ra24}
Ariel Rapaport.
\newblock On self-affine measures associated to strongly irreducible and
  proximal systems.
\newblock {\em Adv. Math.}, 449:Paper No. 109734, 116, 2024.

\end{thebibliography}

\end{document}